\documentclass[12pt,3p,times]{elsarticle}

\usepackage{amsmath,amsfonts,amssymb}
\usepackage{amsthm}
\usepackage{bm}
\usepackage{graphicx,epstopdf}
\usepackage[position=top,labelfont=normalfont,textfont=normalfont,singlelinecheck=off,justification=raggedright]{subfig}
\usepackage{floatrow}
\usepackage[dvipsnames]{xcolor}
\usepackage{setspace}
\usepackage{enumerate,etaremune}
\usepackage{booktabs}
\usepackage{tabularx,multirow}
\usepackage{framed}
\usepackage{hhline}
\usepackage{url}
\usepackage[unicode,bookmarks=false]{hyperref}
\hypersetup{
  colorlinks,
  citecolor=Blue,linkcolor=Blue,urlcolor=Blue
}
\usepackage[toc,page]{appendix}
\usepackage[T1]{fontenc}
\usepackage{cancel}

\usepackage{lineno}

\usepackage[textsize=footnotesize,linecolor=red,backgroundcolor=red!25,bordercolor=red]
  {todonotes}

\usepackage{tikz}
\usetikzlibrary{shapes.geometric}
\usetikzlibrary{shapes,arrows}
\usetikzlibrary{plotmarks}

\usetikzlibrary{calc}

\usepackage{algorithm}
\usepackage{algorithmicx,algpseudocode}
\usepackage{cases}

\newcommand{\eg}{{\it e.g.}}
\newcommand{\ie}{{\it i.e.}}

\newcommand{\etal}{{\it et al.}}
\newcommand{\tensor}[1]{\bm{#1}}

\newcommand{\dd}{\mathrm{d}}
\newcommand{\pd}{\partial}

\newcommand{\rn}[1]{\uppercase\expandafter{\romannumeral #1\relax}}

\DeclareMathOperator{\grad}{\nabla}
\DeclareMathOperator{\diver}{\nabla\cdot}

\DeclareMathOperator{\dyadic}{\otimes}
\newsavebox{\dotbox}

\theoremstyle{remark}
\newtheorem{remark}{Remark}

\newcommand{\revised}[1]{{\color{black} #1}}
\newcommand{\revisedsecond}[1]{{\color{black} #1}}

\newcolumntype{L}[1]{>{\raggedright\let\newline\\arraybackslash\hspace{0pt}}m{#1}}
\newcolumntype{C}[1]{>{\centering\let\newline\\arraybackslash\hspace{0pt}}m{#1}}
\newcolumntype{R}[1]{>{\raggedleft\let\newline\\arraybackslash\hspace{0pt}}m{#1}}

\allowdisplaybreaks

\biboptions{sort&compress,square,comma,numbers}
\AtBeginDocument{\hypersetup{citecolor=MidnightBlue,linkcolor=MidnightBlue,urlcolor=MidnightBlue}}

\usepackage{etoolbox}
\makeatletter
\patchcmd{\ps@pprintTitle}
{Preprint submitted to}
{~}
{}{}
\makeatother

\begin{document}

\begin{frontmatter}

\title{Circumventing volumetric locking in explicit material point methods: A simple, efficient, and general approach}

\author[KAIST]{Yidong Zhao}
\author[UCLA]{Chenfanfu Jiang}
\author[KAIST]{Jinhyun Choo\corref{corr}}
\ead{jinhyun.choo@kaist.ac.kr}

\cortext[corr]{Corresponding Author}

\address[KAIST]{Department of Civil and Environmental Engineering, KAIST, South Korea}
\address[UCLA]{Department of Mathematics, University of California, Los Angeles, United States}

\journal{~}

\begin{abstract}
The material point method (MPM) is frequently used to simulate large deformations of nearly incompressible materials such as water, rubber, and undrained porous media.
However, MPM solutions to nearly incompressible materials are susceptible to volumetric locking, that is, overly stiff behavior with erroneous strain and stress fields.
While several approaches have been devised to mitigate volumetric locking in the MPM, they require significant modifications of the existing MPM machinery, often tailored to certain basis functions or material types.
In this work, we propose a locking-mitigation approach featuring an unprecedented combination of simplicity, efficacy, and generality for a family of explicit MPM formulations. 
The approach combines the assumed deformation gradient ($\bar{\bm{F}}$) method with a volume-averaging operation built on the standard particle--grid transfer scheme in the MPM.
Upon explicit time integration, this combination yields a new and simple algorithm for updating the deformation gradient, preserving all other MPM procedures. 
The proposed approach is thus easy to implement, low-cost, and compatible with the existing machinery in the MPM.
Through various types of nearly incompressible problems in solid and fluid mechanics, we verify that the proposed approach efficiently circumvents volumetric locking in the explicit MPM, regardless of the basis functions and material types.
\end{abstract}


\end{frontmatter}


\section{Introduction}

The material point method (MPM)~\cite{sulsky1994particle} is a hybrid Lagrangian--Eulerian numerical technique for continuum mechanics simulation, whereby physical quantities are traced via material points, or ``particles,'' and the governing equation (in its weak form) is solved in a background grid. 
The particles and the grid exchange their information through projection operations which use basis (interpolation) functions associated with the grid.
Remarkably, the MPM shares many features with the finite element method (FEM)---the most popular numerical method in solid mechanics---while it is rooted in the fluid implicit particle (FLIP) method~\cite{brackbill1986flip}---a particle-in-cell (PIC) method for fluid dynamics.
For this reason, the MPM has been commonly used for simulating large deformation in a wide variety of solids and fluids alike (\eg~\cite{gaume2018dynamic,fern2019material,zhao2020stabilized,li2021three}).

When modeling nearly incompressible materials (\eg~water, rubber, and undrained porous media), MPM solutions are susceptible to volumetric locking, that is, overly stiff behavior with erroneous strain and stress fields.
Volumetric locking commonly arises in the FEM and its related methods when integration points pose excessive incompressibility constraints on the calculation of element stiffness.
Unfortunately, the MPM is inherently vulnerable to volumetric locking, because it typically uses a large number of integration (material) points per element.
In the FEM literature, a number of approaches have been proposed for mitigating locking, relaxing the incompressibility constraints on element kinematics in different ways.
Among them, only those compatible with finite deformation kinematics as well as movable integration points may be adapted to the MPM.

Locking-mitigation approaches that have been adapted to the MPM can be categorized into the following four types.
The first type is mixed multi-field formulations.
For example, Love and Sulsky~\cite{love2006energy} and Mast \etal~\cite{mast2012mitigating} have used a three-field MPM formulation based on the Hu--Washizu variational principle, and Iaconeta \etal~\cite{iaconeta2019stabilized} have presented a two-field formulation with stabilization.
The second type is operator splitting algorithms explored by Zhang \etal~\cite{zhang2017incompressible} and Kularathna and Soga~\cite{kularathna2017implicit}, which are built on Chorin's projection method for incompressible fluid dynamics~\cite{chorin1968numerical}.
The third type is the assumed deformation gradient ($\bar{\bm{F}}$) method~\cite{de1996design}. 
Coombs and coworkers~\cite{coombs2018overcoming,wang2021efficient} have developed implicit $\bar{\bm{F}}$ MPM formulations for quasi-static solid mechanics,
and Moutsanidis \etal~\cite{moutsanidis2020treatment} have proposed a different way to calculate $\bar{\bm{F}}$ in explicit particle methods.
The fourth one is the nonlinear $\bar{\bm{B}}$ method~\cite{simo1985variational}, which is a large-strain generalization of the $\bar{\bm{B}}$ method originally proposed by Hughes~\cite{hughes1980generalization} for the FEM at small strain, extended to the MPM recently by Bisht \etal~\cite{bisht2021simulating,bisht2021material}.
It is noted that the $\bar{\bm{F}}$ and nonlinear $\bar{\bm{B}}$ methods are highly related in that both methods rely on reduced integration of the volumetric part of a deformation measure (but their specific procedures are not the same).
Very recently, Telikicherla and Moutsanidis~\cite{telikicherla2022treatment} have proposed a projection technique for reduced integration in the MPM with high-order basis functions, whereby an additional background grid with low-order basis functions is introduced.

Nevertheless, the aforementioned locking-mitigation approaches require significant modifications of the existing MPM machinery, often tailored to certain basis functions or material types.
The mixed formulations and operator splitting algorithms require one to change the standard governing equations and time-stepping scheme, respectively, demanding significant costs for implementation and utilization. 
Extension of these methods to multiphysical problems (\eg~coupled deformation and flow) is also a challenging endeavor.
Regarding the $\bar{\bm{F}}$ and $\bar{\bm{B}}$ methods, their current MPM versions are either restricted to basis functions that are not associated with adjacent elements (undesirable due to cell-crossing errors) or require a non-trivial modification of basis functions if they are associated with adjacent elements.
For example, an additional basis function dedicated to element-wise averaging is necessary for the $\bar{\bm{F}}$ MPM formulation of Coombs~\etal~\cite{coombs2018overcoming}, which is built on the generalized interpolation material point (GIMP) method.
As such, it is not straightforward to apply the existing $\bar{\bm{F}}$ or $\bar{\bm{B}}$ method to MPM formulations with different MPM basis functions such as B-splines~\cite{steffen2008analysis,gan2018enhancement}. 

In this work, we present a new approach that is unprecedentedly simple, efficient, and general for circumventing volumetric locking in a family of standard explicit MPM formulations.
The key idea is to calculate the assumed deformation gradient, $\bar{\bm{F}}$, using the standard particle--grid transfer scheme in the MPM, instead of the element-wise averaging operation or multiple background grids used in the existing $\bar{\bm{F}}$ methods for MPM.
Combining this idea with the standard explicit time discretization in the MPM, we arrive at a new and simple algorithm for updating the deformation gradient.
The new algorithm neither changes any other parts of the existing MPM machinery nor introduces any additional parameters.
Therefore, this approach can be utilized in a highly straightforward and efficient manner, regardless of the MPM basis functions and material types.
We implement the proposed approach with two types of MPM basis functions, namely, GIMP's basis functions and B-splines, and verify it with various types of nearly incompressible problems arising in solid and fluid mechanics.

\section{Material point method formulation}
\label{section:mpm-formulation}
  
This section recapitulates the standard MPM formulation for a continuum body undergoing large deformation.
For more details of the formulation, the reader is referred to~\cite{jiang2016material,zhang2016material,de2020material}.
It is noted that while there exist a number of more advanced MPM formulations (\eg~\cite{jiang2015affine,fu2017polynomial,hammerquist2017new,jiang2017angular}), their stress update procedure---the part to which our locking-mitigation approach will be applied---is identical to that in the standard MPM formulation. As such, here we shall focus on the basic and commonly used MPM formulation.

\subsection{Problem statement}
Consider a continuum body whose current configuration is denoted by $\Omega \in \mathbb{R}^{\dim}$, where ``$\dim$'' refers to the spatial dimension. 
The boundary of $\Omega$ is denoted by $\pd \Omega$, and it is decomposed into the displacement (Dirichlet) boundary $\pd_{u}\Omega$ and the traction (Neumann) boundary $\pd_{t}\Omega$ such that $\pd_{u}\Omega \cap \pd_{t}\Omega = \emptyset$ and $\overline{\pd_{u}\Omega \cup \pd_{t}\Omega} = \pd\Omega$.
The time domain is denoted by $\mathcal{T}:=(0,T]$ with $T > 0$.

Finite deformation theory should be used to accurately describe nonlinear kinematics in large deformation.
Let us denote by $\tensor{X}$ and $\tensor{x}$ the position vectors of a material point in the reference and current configurations, respectively.
The displacement vector of the material point is then defined as $\tensor{u} := \tensor{x} - \tensor{X}$.
The velocity and acceleration vectors are given by $\tensor{v} := \dot{\bm{u}}$ and $\tensor{a} := \dot{\bm{v}} = \ddot{\bm{u}}$, where the dot denotes the material time derivative.
The deformation gradient is defined as
\begin{equation}
  \tensor{F} := \dfrac{\partial \tensor{x}}{\partial \tensor{X}} = \tensor{1} + \dfrac{\partial \tensor{u}}{\partial \tensor{X}}\,,
\end{equation}
where $\tensor{1}$ is the second-order identity tensor.
The Jacobian is defined as
\begin{equation}
  J := \det{(\tensor{F})} = \dfrac{\dd v}{\dd V},
\end{equation}
where $\dd V$ and $\dd v$ are the differential volumes in the reference and current configurations, respectively.

The balance of linear momentum provides the governing equation.
Since the standard MPM is built on the updated Lagrangian approach, we write the momentum balance equation in the current configuration as
\begin{align}
  \diver\tensor{\sigma}(\bm{F}) + \rho\tensor{g} = \rho\dot{\tensor{v}} \quad
  &\text{in} \:\: \Omega \times \mathcal{T},
  \label{eq:strong-form}
\end{align}
where $\diver\,(\circ)$ denotes the divergence operator defined in the current configuration, 
$\tensor{\sigma}$ is the Cauchy stress tensor,
$\rho$ is the mass density,
and $\tensor{g}$ is the gravitational acceleration vector.
To close the equation, a constitutive relation between $\bm{\sigma}$ and $\bm{F}$ should be introduced.
In this work, we will consider a range of commonly used constitutive relations to demonstrate the generality of the proposed method. 
For brevity, we omit the details of these constitutive relations, referring to textbooks on this subject (\eg~\cite{holzapfel2002nonlinear,de2011computational,borja2013plasticity}).

The initial--boundary value problem of interest can be stated as follows:
Find $\bm{u}$ that satisfies Eq.~\eqref{eq:strong-form}, subject to the initial condition $\tensor{u} = \tensor{u}_{0}$ and boundary conditions 
\begin{align}
    \tensor{u} = \hat{\tensor{u}} \quad
    &\text{on}\:\: \partial_{u}\Omega \times \mathcal{T}, \\
    \tensor{n}\cdot\tensor{\sigma} = \hat{\tensor{t}} \quad
    &\text{on}\:\: \partial_{t}\Omega\times \mathcal{T},
\end{align}
where $\hat{\tensor{u}}$ and $\hat{\tensor{t}}$ are the boundary displacement and traction, respectively, and $\tensor{n}$ is the unit outward normal vector in the current configuration.

Through the standard weighted residual procedure, the variational form of the governing equation can be written as
\begin{equation}
  \int_{\Omega} \tensor{\eta} \cdot \rho \dot{\tensor{v}} \, \dd V
  =
  -\int_{\Omega} \grad^{\mathrm{s}} \tensor{\eta} : 
  \tensor{\sigma}(\bm{F}) \, \dd V
  +
  \int_{\Omega} \tensor{\eta} \cdot \rho \tensor{g} \, \dd V
  +
  \int_{\partial_{t} \Omega} \tensor{\eta} \cdot \hat{\tensor{t}} \, \dd A,
  \label{eq:variational-equation}
\end{equation}
where $\tensor{\eta}$ denotes the variation of the displacement field, and $\grad^{\mathrm{s}}$ is the symmetric gradient operator defined in the current configuration.

\subsection{Material point method discretization}
\begin{figure}[htbp]
  \centering
  \includegraphics[width=0.9\textwidth]{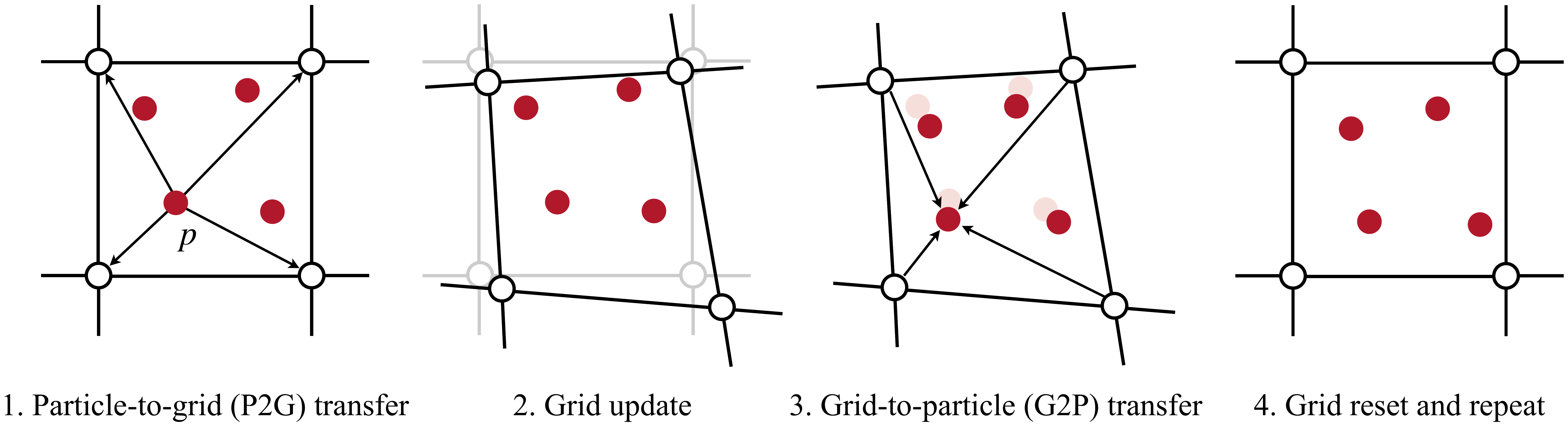}
  \caption{MPM update procedure.}
  \label{fig:mpm-procedure}
\end{figure}

For MPM discretization of the problem, we introduce a set of particles (material points) filling in the domain and a background grid that accommodates the particles.
We then update the solution from the previous time step ($t^{n}$) to the next time step ($t^{n+1}$) through the procedure illustrated in Fig.~\ref{fig:mpm-procedure}.
The MPM procedure is recapitulated in the following. 
Hereafter, we shall use subscript $(\circ)_{p}$ to denote quantities related to particles and use subscript $(\circ)_{i}$ to denote quantities related to nodes.
We shall also use superscripts $(\circ)^{n}$ and $(\circ)^{n+1}$ to denote quantities at $t^{n}$ and $t^{n+1}$, respectively.

\paragraph{Particle-to-grid (P2G) transfer} 
First, we map the mass and momentum of each particle to nodes in the background grid.
This process is called the particle-to-grid (P2G) transfer.
The particle mass and momentum are transferred to the nodes as
\begin{align}
  m_i &= \sum_{p} w_{ip} m_p, \label{eq:P2G-mass}\\
  m_i\bm{v}_i^{n} &= \sum_{p} w_{ip} m_p \bm{v}_p^{n}, \label{eq:P2G-momentum}
\end{align}
where $m_{i}$ and $m_p$ are the nodal and particle masses, respectively, and $\bm{v}_{i}$ and $\bm{v}_{p}$ are the particle and nodal velocity vectors, respectively. 
Also,  $w_{ip} := w_i (\tensor{x}_p)$ is the basis function for interpolating values at node $i$ to the position of particle $p$, and $\sum_{p}$ is the summation over particles supported by the basis function associated with node $i$.
In the MPM literature, a few types of basis functions have been employed.
Among them, here we consider two popular choices, namely, (i) the basis functions of the generalized interpolation material point (GIMP) method~\cite{bardenhagen2004generalized}, and (ii) B-splines~\cite{steffen2008analysis}.
Note that these basis functions are free of the cell-crossing error problem which may arise when the original MPM basis functions are used.

\paragraph{Grid update} 
Following the P2G transfer, we update the velocity vector at each node through the governing equation~\eqref{eq:variational-equation}.
As standard, we use the explicit Euler method to integrate the governing equation in time and update the nodal velocity as
\begin{equation}
  m_{i}\bm{v}_{i}^{n+1} = m_{i}\bm{v}_{i}^{n} + \Delta{t}(\bm{f}^{\mathrm{int}}_i + \bm{f}^{\mathrm{ext}}_i).
  \label{eq:nodal-update}
\end{equation}
Here, $\Delta{t} := t^{n+1} - t^{n}$ is the time increment, and $\bm{f}^{\mathrm{int}}_i$ and $\bm{f}^{\mathrm{ext}}_i$ are the internal and external force vectors, respectively.
The internal force vector is calculated as
\begin{align}
  \bm{f}^{\mathrm{int}}_i &= -\sum_{p} \grad^{\mathrm{s}} w_{ip} : \tensor{\sigma} (\bm{F}_p)\, V_p^{n},
  \label{eq:internal-force-vector}
\end{align}
where $V_p^{n}$ is the current particle volume.
The external force vector is calculated from the body force and boundary traction vectors.

\paragraph{Grid-to-particle (G2P) transfer} 
After updating the nodal velocity, we map it back to the particles and update the particle velocity.
This procedure is called the grid-to-particle (G2P) transfer.
Two schemes exist for updating the particle velocity: 
(i) the fluid-implicit-particle (FLIP) method~\cite{brackbill1986flip} which transfers the velocity increment, 
and (ii) the particle-in-cell (PIC) method~\cite{harlow1964particle} which transfers the updated velocity itself.
These two schemes can be written as
\begin{align}
  (\tensor{v}_p^{n+1})_{\text{FLIP}} &= \tensor{v}_p^n + \sum_{i} w_{ip} (\tensor{v}^{n+1}_i - \tensor{v}_i^n), \label{eq:FLIP} \\
  (\tensor{v}_p^{n+1})_{\text{PIC}} &= \sum_{i} w_{ip}\tensor{v}^{n+1}_i, \label{eq:PIC}
\end{align}
where $\sum_{i}$ is the summation over nodes supporting particle $p$.
While FLIP is less stable than PIC, it has significantly less numerical damping than PIC.
In general, the FLIP and PIC schemes can be blended as
\begin{equation}
  \tensor{v}_p^{n+1} = \eta (\tensor{v}_p^{n+1})_{\text{FLIP}} + (1 - \eta)(\tensor{v}_p^{n+1})_{\text{PIC}}.
\end{equation}
where $\eta\in[0,1]$ is the blending coefficient. Unless otherwise specified, we shall use $\eta=1$.
Subsequently, we update the deformation gradient, volume, stress, and position of each particle.
The deformation gradient is updated as
\begin{align}
  \bm{F}_{p}^{n+1} = \left(\tensor{1} + \Delta t \sum_{i} \tensor{v}_i^{n+1}\dyadic\grad w_{ip}\right)\cdot\bm{F}_{p}^{n}.
  \label{eq:deformation-gradient-update}
\end{align} 
The particle volume is updated as
\begin{equation}
  V_{p}^{n+1} = J_{p}^{n+1}V_{p}^{0},
  \label{eq:particle-volume-update}
\end{equation}
where $J_{p}^{n+1} := \det{\tensor{F}}_{p}^{n+1}$, and $V_{p}^{0}$ is the initial (reference) particle volume. 
For updating the stress tensor, we use the \emph{relative} deformation gradient---analogous to the incremental strain tensor in small-strain problems---to accommodate history-dependent material behavior.
The relative deformation gradient is defined and calculated as
\begin{equation}
  \Delta \tensor{F}_p
  := {\tensor{F}}_{p}^{n+1} \cdot ({\tensor{F}}_{p}^{n})^{-1}
  = \tensor{1} + \Delta t \sum_{i} \tensor{v}_i^{n+1}\dyadic\grad w_{ip}.
  \label{eq:relative-deformation-gradient}
\end{equation}
By assigning the relative deformation gradient to the specific constitutive relation, the stress tensor is updated.
Lastly, the particle position is updated as
\begin{equation}
  \tensor{x}_p^{n+1} = \tensor{x}_p^n + \Delta t \sum_{i} w_{ip}\tensor{v}^{n+1}_i. 
\end{equation}

\paragraph{Grid reset and repeat}
After updating the state variables and positions of the particles, we reset the background grid, move to the next time step, and repeat the aforementioned procedure.
It would be worthwhile to note that the grid is not reset in total Lagrangian MPM formulations (\eg~\cite{deVaucorbeil2020total,deVaucorbeil2021modelling}).
Here, we focus on the standard MPM formulations built on an updated Lagrangian approach.

\section{Circumventing volumetric locking}
\label{section:locking-free-mpm}

In this section, we formulate a new approach for mitigating volumetric locking in the MPM.
We first adopt the assumed deformation gradient ($\bar{\bm{F}}$) method---originally proposed for overcoming volumetric locking in FEM~\cite{de1996design}---in the context of the MPM.
We then present a new way to calculate the assumed deformation gradient in the MPM, which averages the volumetric part of the deformation gradient (the Jacobian) through the standard particle--grid transfer operations.
Subsequently, we develop a detailed procedure to apply the proposed method in the MPM.

\subsection{Assumed deformation gradient method}
The key idea of the assumed deformation gradient method is to replace the deformation gradient ($\bm{F}$) in the constitutive relation with an assumed deformation gradient ($\bar{\bm{F}}$), where the volumetric part of $\bm{F}$ (\ie~the Jacobian, $J$) is volume-averaged in some manner.
When applied to the current formulation, it replaces $\bm{F}_{p}$ in the constitutive relation by
\begin{equation}
  \bar{\tensor{F}}_p = \left( \dfrac{\bar{J}_p}{J_p} \right)^{1/\dim} \tensor{F}_p,
  \label{eq:F-bar}
\end{equation}
where $\bar{J}_{p}$ denotes the averaged Jacobian, which is subject to fewer volumetric constraints than $J_p$.
It can be seen that $(\bar{J}_p/J_p)^{1/\dim}$ acts as a scaling term for the deformation gradient used for the constitutive update.

In the FEM, for which the $\bar{\bm{F}}$ method was originally proposed, $\bar{J}_p$ can be calculated straightforwardly as the average of $J_{p}$ in each element.
In the MPM, however, such element/cell-wise averaging is not ideal, because the basis functions are often related to adjacent elements to avoid cell-crossing errors (\eg~GIMP's basis functions and B-splines).
As such, when Coombs~\etal~\cite{coombs2018overcoming} applied the $\bar{\bm{F}}$ method to GIMP's basis functions, they had to introduce an additional basis function specialized in element-wise averaging.
Unfortunately, introducing such a special basis function is not only cumbersome but also restricted to a specific MPM scheme.
For example, the formulation in Coombs~\etal~\cite{coombs2018overcoming} is not compatible with other MPM basis functions.

\subsection{Volume averaging}
In this work, we present a new approach to volume-averaging $J_p$ in the MPM, which builds on the existing particle--grid transfer schemes and hence preserves the existing basis functions.
Concretely, it first projects the Jacobians at the particles to the background grid in a volume-averaging manner, in a way similar to the P2G transfer.
Then, through the G2P transfer, the volume-averaged Jacobians at the nodes are mapped back to the particles, so that they can be used for constitutive updates at the individual particles.

To express the approach mathematically, let us define the volume-averaged projection of $J_p$ to the background grid, as
\begin{equation}
  \bar{J}_{i} = \sum_{p} w_{ip} J_p V_p / V_i, \quad V_i := \sum_{p} w_{ip} V_p.
  \label{eq:J_bar_node}
\end{equation}
It can be seen that this projection is more or less the same as the standard P2G transfer, except that the particle volume is considered for volume averaging.
Next, we project back $\bar{J}_{i}$---the volume-averaged Jacobian defined at the grid nodes---to the particles where $\bar{\bm{F}}$ is used for updating the stress tensor. 
For this purpose, we use the standard G2P transfer as 
\begin{equation}
  \bar{J}_p = \sum_{i} w_{ip} \bar{J}_{i}.
\end{equation}
For notational simplicity in the succeeding formulations, we shall express the foregoing operation as the operator $\Pi(\circ)$, say,
\begin{equation}
  \Pi(J_p) := \bar{J}_{p} = \sum_{i} w_{ip} \bar{J}_{i}.
  \label{eq:volume-averaging-operator}
\end{equation}
\smallskip

\begin{remark}
  While the projection-based averaging described above appears similar to that in Telikicherla and Moutsanidis~\cite{telikicherla2022treatment}, there are a couple of important differences.
  First, the projection operation in Telikicherla and Moutsanidis~\cite{telikicherla2022treatment} utilizes additional basis functions whose order is lower than their default basis functions (B-splines) for other parts of the MPM. 
  For this reason, their projection entails an additional background mesh dedicated to the projection.
  Second, their methods apply the projection to the divergences of the velocity and stress fields in the variational equation, whereas we apply the projection to the Jacobian before the stress update.
\end{remark}
\smallskip

\begin{remark}
  The foregoing volume-averaging operation is analogous to the projection scheme utilized for mitigating locking in a different meshfree method in Ortiz-Bernardin \etal~\cite{ortiz2015improved}, where the authors adapted the assumed deformation gradient operation in Broccardo \etal~\cite{broccardo2009assumed} to the meshfree context.
  It is noted, however, that the formulation of Ortiz-Bernardin \etal~\cite{ortiz2015improved} involves a modification of the strain--displacement matrix in addition to the volume-averaging projection. 
  This is different from the approach proposed herein, whereby the volume-averaging operation is used to evaluate $\bar{\bm{F}}$ only and the existing discretization is retained.
\end{remark}
\smallskip

\begin{remark}
  Unlike the existing $\bar{\bm{F}}$ MPM formulations where the volume averaging of $J$ is performed inside individual cells~\cite{coombs2018overcoming,moutsanidis2020treatment,wang2021efficient}, here the volume averaging is done inside the support of the basis functions associated with individual particles.
  This way allows us to accommodate particles that influence multiple cells without any change in the existing basis functions.
  When an implicit integration is used, however, it may be less desirable than introducing a new basis function (as proposed by Coombs \etal~\cite{coombs2018overcoming}), because it would be more onerous to calculate the derivative of $\bar{\bm{F}}$.
  For an explicit integration---dominant in the MPM community---the projection operation must be far simpler than modifying the basis functions.
\end{remark}

\subsection{Stress update procedure}
We now discuss how to update the stress tensor of a material point with the proposed approach.
Consider the stress update stage during an MPM update between $t^{n}$ (previous time step) and $t^{n+1}$ (next time step).
All the quantities of the particle at $t^{n}$, including the assumed deformation gradient, $\bar{\bm{F}}^{n}_{p}$, are known.
However, the quantities at $t^{n+1}$ are unknown, except the deformation gradient, $\bm{F}^{n+1}_{p}$, calculated from Eq.~\eqref{eq:deformation-gradient-update}.

Let us recall that the stress update is based on the relative deformation gradient, Eq.~\eqref{eq:relative-deformation-gradient}.
Since we use the $\bar{\bm{F}}$ method, we define the \emph{relative} $\bar{\bm{F}}$ as
\begin{equation}
  \Delta \bar{\tensor{F}}_p := \bar{\tensor{F}}_{p}^{n+1} \cdot (\bar{\tensor{F}}_{p}^{n})^{-1}.
  \label{eq:relative-assumed-deformation-gradient-v1}
\end{equation}
It is noted that while the (original) deformation gradient at $t^{n+1}$, ${\bm{F}}^{n+1}_{p}$, is given from the updated velocity, the assumed deformation gradient at $t^{n+1}$, $\bar{\tensor{F}}_{p}^{n+1}$, is not given directly.
So we first derive the following expression for $\bar{\tensor{F}}_{p}^{n+1}$:
\begin{align}
  \bar{\tensor{F}}_{p}^{n+1} 
  &= \left( \dfrac{\bar{J}_{p}^{n+1}}{J_{p}^{n+1}} \right)^{1/\dim} \tensor{F}_{p}^{n+1}  \nonumber \\
  &= \left( \dfrac{\bar{J}_{p}^{n+1}}{J_{p}^{n} \Delta J_p} \right)^{1/\dim} \tensor{F}_{p}^{n+1}  \nonumber \\
  &= \left( \dfrac{\bar{J}_{p}^{n+1}  \bar{J}_{p}^{n}}{J_{p}^{n} \bar{J}_{p}^{n} \Delta J_p  } \right)^{1/\dim} \tensor{F}_{p}^{n+1}  \nonumber \\
  &= \left( \dfrac{\bar{J}_{p}^{n+1}}{\bar{J}_{p}^{n} \Delta J_p} \right)^{1/\dim} \Delta \tensor{F}_p  \cdot \tensor{F}_{p}^{n} \left( \dfrac{\bar{J}_{p}^{n}}{J_{p}^{n}} \right)^{1/\dim} \nonumber \\
  &= \left( \dfrac{\bar{J}_{p}^{n+1}}{\bar{J}_{p}^{n} \Delta J_p} \right)^{1/\dim} \Delta \tensor{F}_p \cdot\bar{\tensor{F}}_{p}^{n}
  \label{eq:assumed-deformation-gradient-next-step}
\end{align}
where
\begin{equation}
  \Delta J_p := \det{(\Delta \tensor{F}_p)}, \quad
  \bar{J}_{p}^{n} := \det{(\bar{\tensor{F}}_{p}^{n})}, \quad
  J_p^{n} := \det{(\tensor{F}_{p}^{n})}.
  \label{eq:J-definitions}
\end{equation}

Equation~\eqref{eq:assumed-deformation-gradient-next-step} has replaced $\tensor{F}_{p}^{n+1}$ by the product of $\Delta \tensor{F}_{p}$, which is given in the G2P stage as in Eq.~\eqref{eq:relative-deformation-gradient}, and $\bar{\tensor{F}}_{p}^{n}$, which was used for evaluating the stress tensor at $t^{n}$.
In doing so, it has also replaced $J_{p}^{n+1}$ in the denominator of the scaling term by $\bar{J}_{p}^{n} \Delta J_p$.
Inserting Eq.~\eqref{eq:assumed-deformation-gradient-next-step} into Eq.~\eqref{eq:relative-assumed-deformation-gradient-v1} gives
\begin{align}
  \Delta \bar{\tensor{F}}_p 
  &= \left( \dfrac{\bar{J}_{p}^{n+1}}{\bar{J}_{p}^{n} \Delta J_p} \right)^{1/\dim} \Delta \tensor{F}_p \cdot \bar{\tensor{F}}_{p}^{n} \cdot (\bar{\tensor{F}}_{p}^{n})^{-1} \nonumber \\
  &= \left( \dfrac{\bar{J}_{p}^{n+1}}{\bar{J}_{p}^{n} \Delta J_p}\right)^{1/\dim} \Delta \tensor{F}_p.
  \label{eq:relative-assumed-deformation-gradient-v2}
\end{align} 
Comparing the above equation with Eq.~\eqref{eq:F-bar}, one can see that the above equation is an incremental version of the $\bar{\tensor{F}}$ method.
The final task is to evaluate $\bar{J}_{p}^{n+1}$ as a volume-averaged version of the term in the denominator of the scaling term.
To this end, we apply the volume-averaging operator proposed earlier, $\Pi(\circ)$, as
\begin{equation}
  \bar{J}_{p}^{n+1} = \Pi (\bar{J}_{p}^{n}\Delta J_p).
  \label{eq:J-approximation-explicit}
\end{equation}
Substituting Eq.~\eqref{eq:J-approximation-explicit} into Eq.~\eqref{eq:relative-assumed-deformation-gradient-v2} gives
\begin{equation}
  \Delta \bar{\tensor{F}}_p 
  = \left( \dfrac{\Pi(\bar{J}_{p}^{n} \Delta J_p)}{\bar{J}_{p}^{n} \Delta J_p} \right)^{1/\dim} 
  \Delta \tensor{F}_p.
  \label{eq:relative_bar_deformation_gradient_explicit}
\end{equation}
  
Algorithm~\ref{algo:explicit_Fbar} presents a detailed stress-update procedure in which the proposed locking-mitigation approach is applied.
We emphasize that the proposed approach only modifies the stress update during the G2P stage, preserving all the other aspects described in Section~\ref{section:mpm-formulation}.
For example, the particle volume update, Eq.~\eqref{eq:particle-volume-update}, remains unchanged.

\begin{algorithm}[h!]
  \caption{Stress update procedure with the proposed locking-mitigation approach}
  \begin{enumerate}
    \item Calculate the (original) relative deformation gradient: 
    \begin{equation*}
      \Delta \tensor{F}_p = \tensor{1} + \Delta t \sum_{i} \tensor{v}_i^{n+1}\dyadic\grad w_{ip}.
    \end{equation*}
    \item Compute $\Delta J_p = \det{(\Delta \tensor{F}_p)}$ and $\bar{J}_p^{n} = \det{(\bar{\boldsymbol{F}}_p^{n})}$.
    \item Calculate volume-averaged Jacobians at nodes: 
    \begin{equation*}
        \bar{J}_{i}^{n+1}  = \sum_{p} w_{ip} V_p^{n} (\bar{J}_p^{n}\Delta{J}_{p}) / V_i^{n}, \quad V_i^{n} = \sum_{p} w_{ip} V_p^{n}.
    \end{equation*}
    \item Compute the relative assumed deformation gradient:
    \begin{align*}
      \Delta \bar{\tensor{F}}_p
      &= \left( \dfrac{\Pi(\bar{J}_p^{n}\Delta J_p)}{\bar{J}_p^{n}\Delta J_p} \right)^{1/\dim} \Delta \tensor{F}_p \\
      &= \left( \dfrac{\sum_{i} w_{ip} \bar{J}_{i}^{n+1}}{\bar{J}_p^{n}\Delta J_p} \right)^{1/\dim} \Delta \tensor{F}_p.
    \end{align*}
    \item Update the particle stress with the relative assumed deformation gradient, $\Delta \bar{\tensor{F}}_p$.
  \end{enumerate}
  \label{algo:explicit_Fbar}
\end{algorithm}

\section{Numerical examples}
\label{section:examples}

In this section, we verify and demonstrate the performance of the proposed locking-mitigation approach through four numerical examples involving various types of nearly incompressible materials.
The first example is Cook's membrane~\cite{cook1974improved}, which is a popular benchmark problem for incompressible elasticity.
The second example is the problem of a strip footing on an incompressible elastoplastic solid, for which an analytical solution (the Prandtl solution) is available for the bearing capacity.
The third example is the dam break problem in Mast \etal~\cite{mast2012mitigating}, where a nearly incompressible fluid (water) is allowed to flow freely.
The fourth and last example is a 3D landslide problem in which undrained clay---an elastoplastic and incompressible solid---collapses.
Except for the last 3D example, plane-strain conditions are considered.
Gravity is neglected in the first and second examples.

To confirm that the proposed approach works well regardless of the MPM basis functions, we simulate each example with two different basis functions, namely, GIMP's basis functions and B-splines.
The specific GIMP scheme used herein is uGIMP~\cite{wallstedt2008evaluation} in which the influence domains of individual particles are fixed.
The particular B-splines used in the following results are quadratic B-splines.
While not presented for brevity, we have also found that the proposed method manifests similar performance when cubic B-splines are used.
The MPM results in this section are produced using the \texttt{Taichi} library~\cite{hu2019taichi}.

\subsection{Cook's membrane}
\begin{figure}[htbp]
  \centering
  \includegraphics[width=0.4\textwidth]{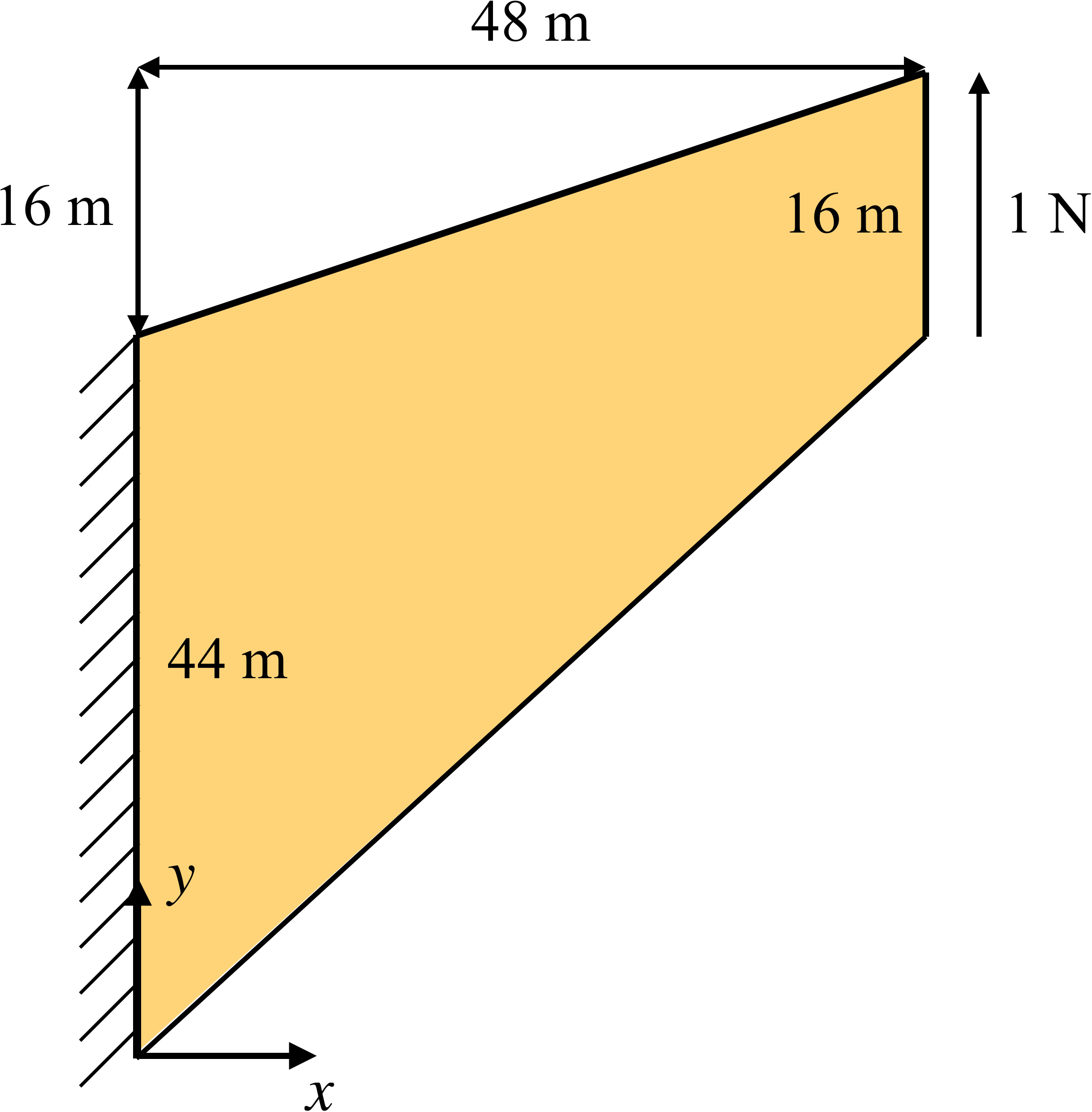}
  \caption{Cook's membrane: problem geometry and boundary conditions.}
  \label{fig:cooks_membrane_setup}
\end{figure}

To verify our formulation in a simple setting, we first simulate Cook's membrane problem~\cite{cook1974improved}, which has widely been used as a benchmark problem for incompressible elasticity.
\revisedsecond{It is noted that while the problem does not have an analytical solution, its benchmark solutions have been produced by several types of locking-free numerical methods (\eg~\cite{rodriguez2016particle,iaconeta2019stabilized,bisht2021simulating}).}
As shown in Fig.~\ref{fig:cooks_membrane_setup}, this problem considers a trapezoidal membrane subjected to a distributed shear traction on its right side and clamped on its left side.
To compare our $\bar{\tensor{F}}$ MPM results with those of $\bar{\tensor{B}}$ MPM results in Bisht \etal~\cite{bisht2021simulating}, we set the shear load and material parameters identical to those in the reference paper.
The load is set as 1 N.
The membrane is a Neo-Hookean solid with Young's modulus of $E = 70$ Pa and Poisson's ratio of $\nu = 0.499$ (nearly incompressible).

To examine the convergence of the numerical solutions, we mainly use two levels of spatial discretization: (i) a coarse discretization that uses 5,776 material points with a background grid comprised of 1-m-long square elements, and (ii) a fine discretization that uses 23,072 material points with a background grid comprised of 0.5-m-long square elements.
To simulate this quasi-static problem with the current explicit dynamic formulation, we adopt the local damping method in Al-Kafaji~\cite{alkafaji2013formulation}, with the same local damping factor used in Bisht \etal~\cite{bisht2021simulating}.
We calculate the time increment as $\Delta t = 0.3 (h/c)$, where $h$ is the element size and $c$ is the P-wave velocity.
This calculation gives $\Delta t = 2.773\times 10^{-3}$ s for the coarse grid and $\Delta t = 1.386\times 10^{-3}$ s for the fine grid.
To ensure numerical stability, we set the FLIP/PIC blending ratio as $\eta = 0.85$.

\begin{figure}[htbp]
  \centering
  \subfloat[Coarse discretization (5,776 material points)]{\includegraphics[width=1.0\textwidth]{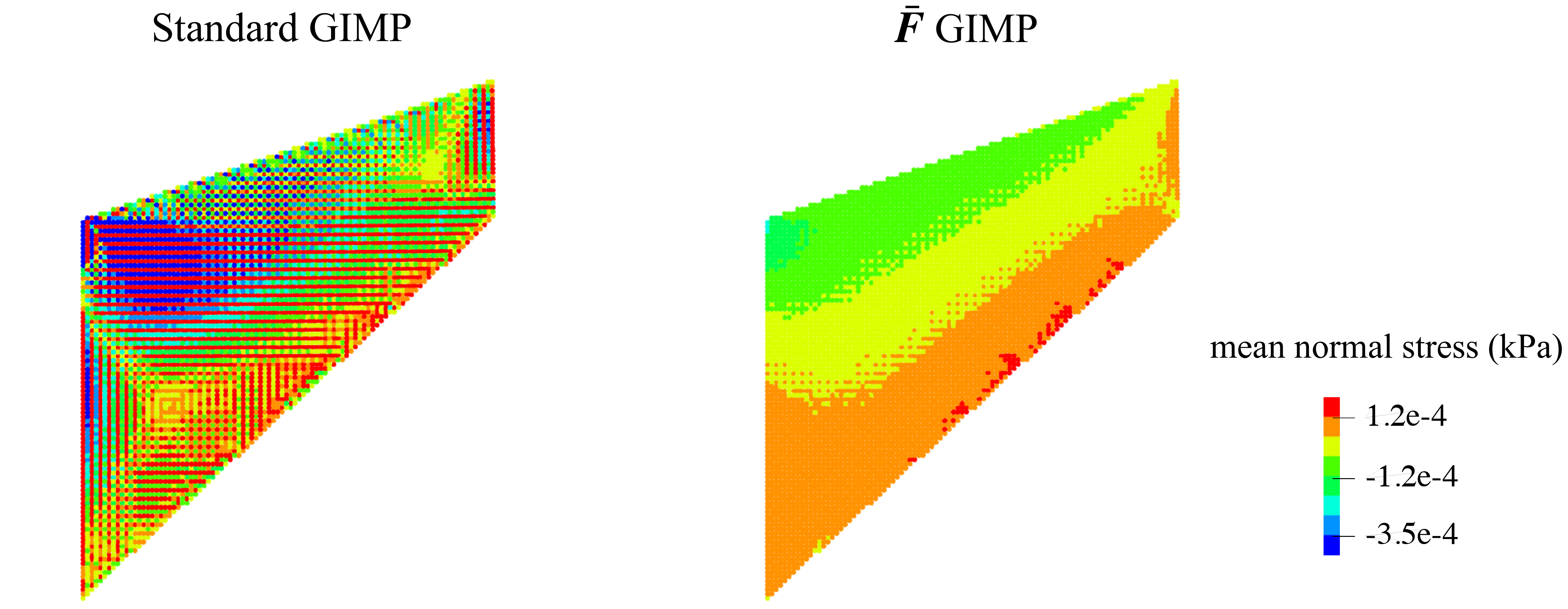}} \newline
  \subfloat[Fine discretization (23,072 material points)]{\includegraphics[width=1.0\textwidth]{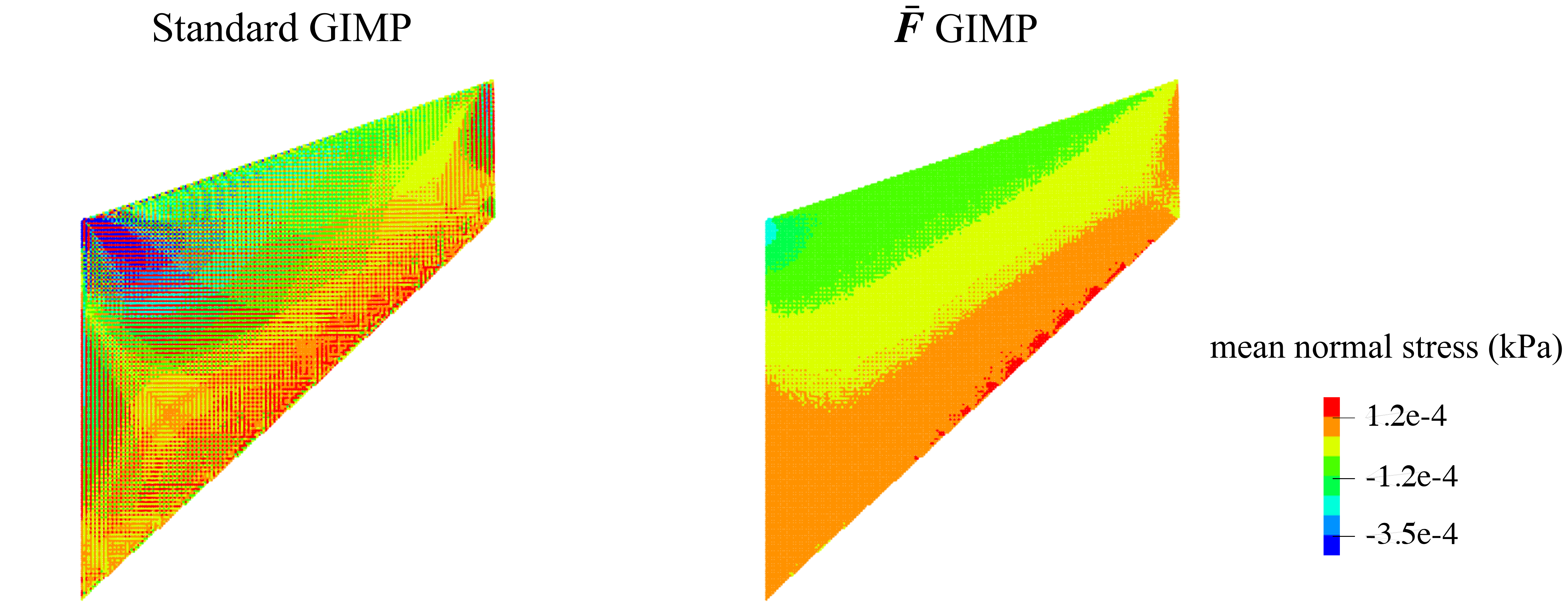}}
  \caption{Cook's membrane: mean normal stress fields in the standard and $\bar{\tensor{F}}$ MPM solutions, obtained with GIMP basis functions.}
  \label{fig:cooks_membrane_GIMP}
\end{figure}

Figure~\ref{fig:cooks_membrane_GIMP} presents the mean normal stress fields obtained by the standard and $\bar{\tensor{F}}$ MPM formulations with GIMP.
As can be seen, the standard MPM solutions are plagued by severe non-physical stress oscillations, which are not remedied by spatial refinement.
Such oscillations have been commonly observed in numerical solutions affected by volumetric locking.
Meanwhile, the $\bar{\tensor{F}}$ MPM solutions are free of non-physical oscillations in the stress fields.
This difference indicates that the proposed formulation does not suffer from volumetric locking.
It is noted that although the $\bar{\tensor{F}}$ MPM solution still has minor oscillations, these remaining oscillations are unrelated to volumetric locking. This will be confirmed later through a comparison between our $\bar{\tensor{F}}$ MPM solution with a nonlinear $\bar{\tensor{B}}$ MPM solution to the same problem~\cite{bisht2021simulating}.

\begin{figure}[htbp]
  \centering
  \subfloat[Coarse discretization (5,776 material points)]{\includegraphics[width=1.0\textwidth]{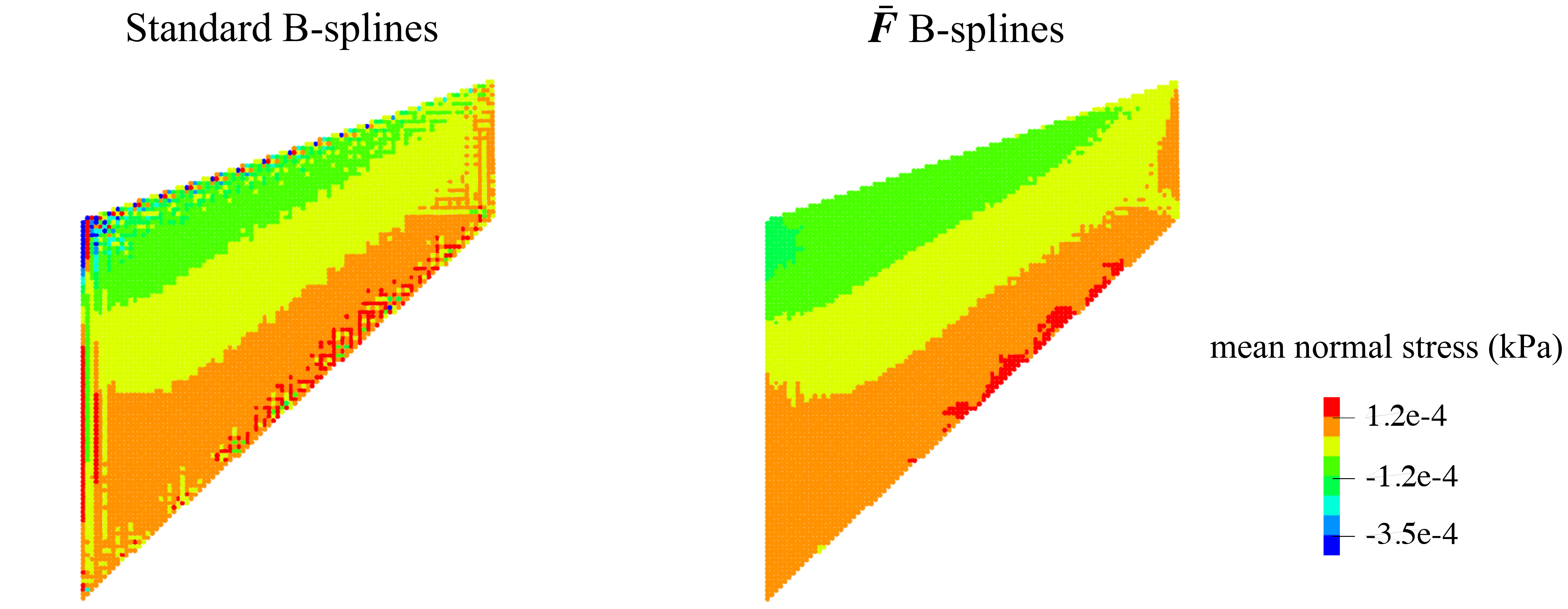}} \newline
  \subfloat[Fine discretization (23,072 material points)]{\includegraphics[width=1.0\textwidth]{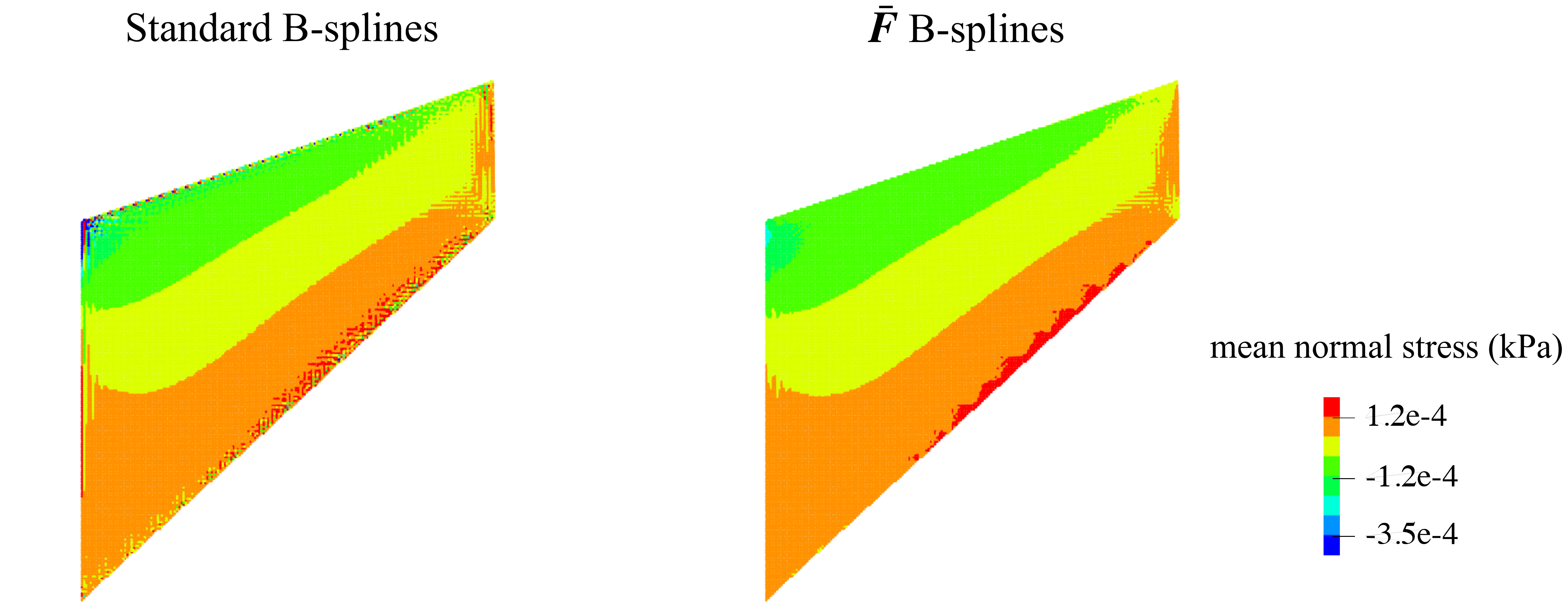}}
  \caption{Cook's membrane: mean normal stress fields in the standard and $\bar{\tensor{F}}$ MPM solutions, obtained with B-splines basis functions.}
  \label{fig:cooks_membrane_Bspline}
\end{figure}

Figure~\ref{fig:cooks_membrane_Bspline} shows how the mean normal stress fields become different when the basis functions are changed to B-splines.
One can see that stress fields in the B-splines MPM are less oscillatory than those in GIMP, especially when the discretization is fine.
This implies that quadratic B-splines are less constrained than GIMP's basis function by the same number of material (integration) points.
Still, however, the standard MPM solutions show undesirable oscillations, because the MPM uses quite a large number of material points per element.
Meanwhile, the same $\bar{\tensor{F}}$ MPM formulation continues to work well notwithstanding the change in the basis functions.

\begin{figure}[htbp]
  \centering
  \includegraphics[width=1.0\textwidth]{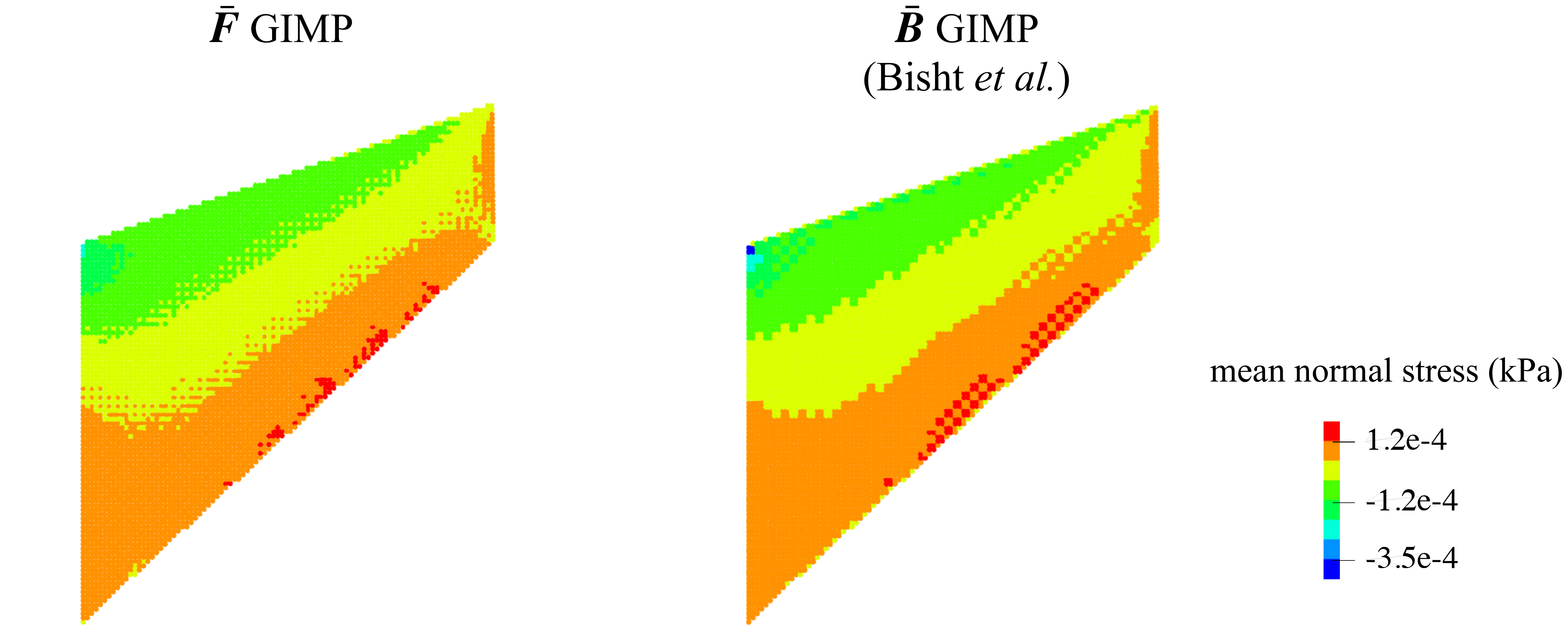}
  \caption{\revised{Cook's membrane: comparison of our $\bar{\tensor{F}}$ GIMP solution with the nonlinear $\bar{\tensor{B}}$ GIMP solution in Bisht \etal~\cite{bisht2021simulating}. Both solutions are produced with 5,776 material points.}}
  \label{fig:cooks_membrane_verification}
\end{figure}

\begin{figure}[htbp]
  \centering
  \subfloat[GIMP]{\includegraphics[width=0.45\textwidth]{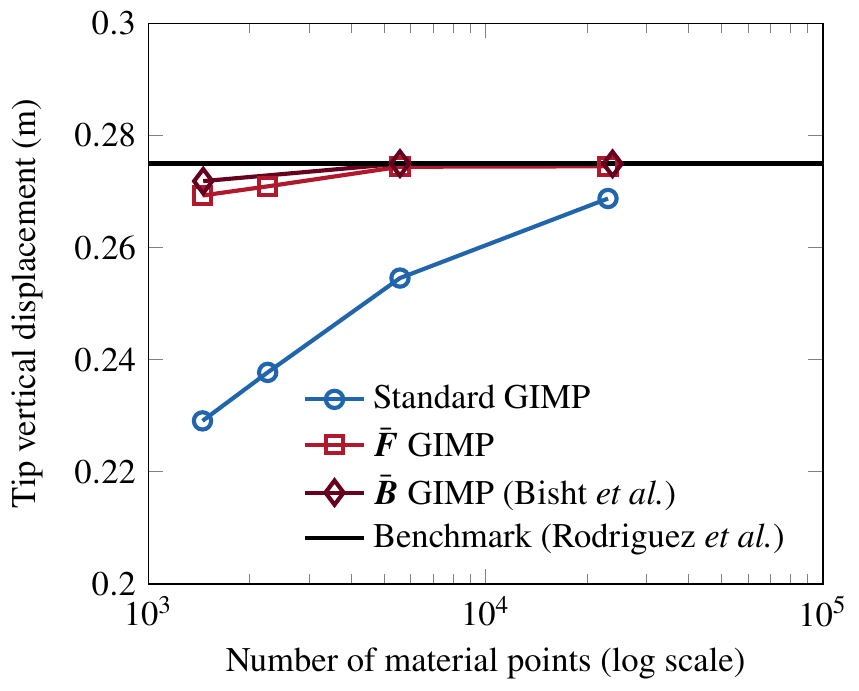}\label{fig:cooks_membrane_tip_displacement_gimp}}\hspace{0.5em}
  \subfloat[B-splines]{\includegraphics[width=0.45\textwidth]{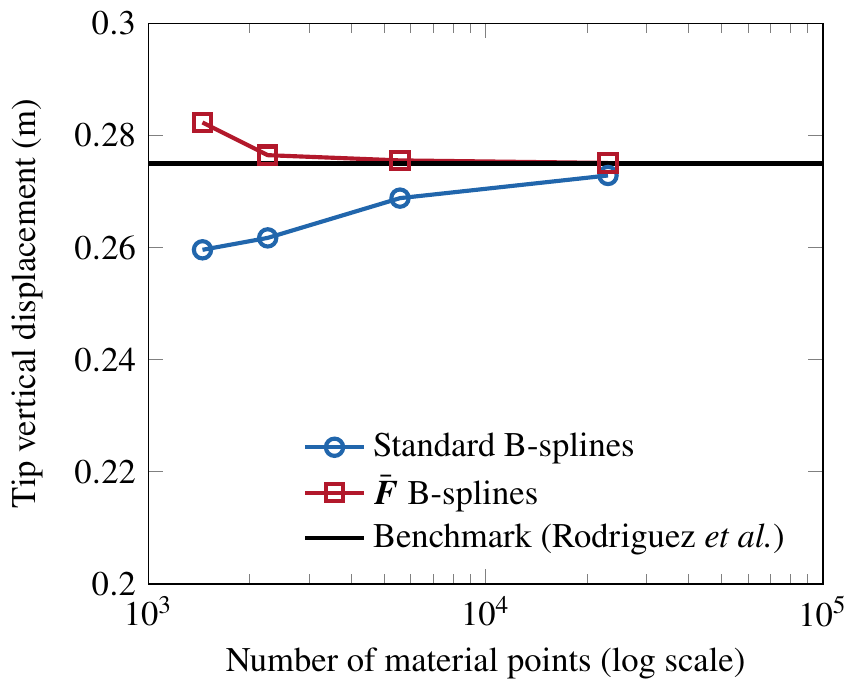}\label{fig:cooks_membrane_tip_displacement_bsplines}}
  \caption{\revised{Cook's membrane: tip vertical displacements from our standard and $\bar{\tensor{F}}$ MPM solutions, along with the nonlinear $\bar{\bm{B}}$ GIMP solution in Bisht \etal~\cite{bisht2021simulating}, calculated with a different number of material points. The benchmark solution is obtained by Rodriguez \etal~\cite{rodriguez2016particle} using an extremely fine discretization.}}
  \label{fig:cooks_membrane_tip_displacement}
\end{figure}

For further verification, we compare our GIMP solution with the reference solution from Bisht \etal~\cite{bisht2021simulating}---obtained with a nonlinear $\bar{\tensor{B}}$ method specialized to GIMP's basis functions---in terms of the mean normal stress field and the vertical displacement at the tip (upper right corner), respectively, in Figs.~\ref{fig:cooks_membrane_verification} and~\ref{fig:cooks_membrane_tip_displacement_gimp}.
It would be worthwhile to note that the nonlinear $\bar{\tensor{B}}$ MPM solution itself has also been verified with other results in the literature, see Bisht \etal~\cite{bisht2021simulating} for more details.
From Fig.~\ref{fig:cooks_membrane_verification}, we can see that the mean normal stress fields in the two solutions are very similar. 
Also, the fact that the $\bar{\tensor{F}}$ and $\bar{\tensor{B}}$ solutions exhibit similar oscillations indicates that these oscillations are unrelated to volumetric locking.
Figure~\ref{fig:cooks_membrane_tip_displacement_gimp} also shows that the tip displacements in the $\bar{\tensor{F}}$ and $\bar{\tensor{B}}$ solutions are very \revised{similar whenever the number of material points is sufficient.
Importantly, upon spatial refinement, both solutions converge gradually to the benchmark solution obtained with an extremely fine discretization by Rodriguez \etal~\cite{rodriguez2016particle}.} 
\revised{Meanwhile, the tip displacements} in the standard MPM solutions are noticeably lower than the locking-free solutions.
\revised{For completeness, in Fig.~\ref{fig:cooks_membrane_tip_displacement_bsplines} we also show the tip displacements in the solutions obtained with B-splines basis functions.
As can be seen, the $\bar{\tensor{F}}$ MPM solutions with B-splines also converge to the benchmark solution upon spatial refinement.}
Taken together, it has been confirmed that the proposed method performs similarly to the nonlinear $\bar{\tensor{B}}$ method.

\subsection{Strip footing}

\begin{figure}[htbp]
  \centering
  \includegraphics[width=0.45\textwidth]{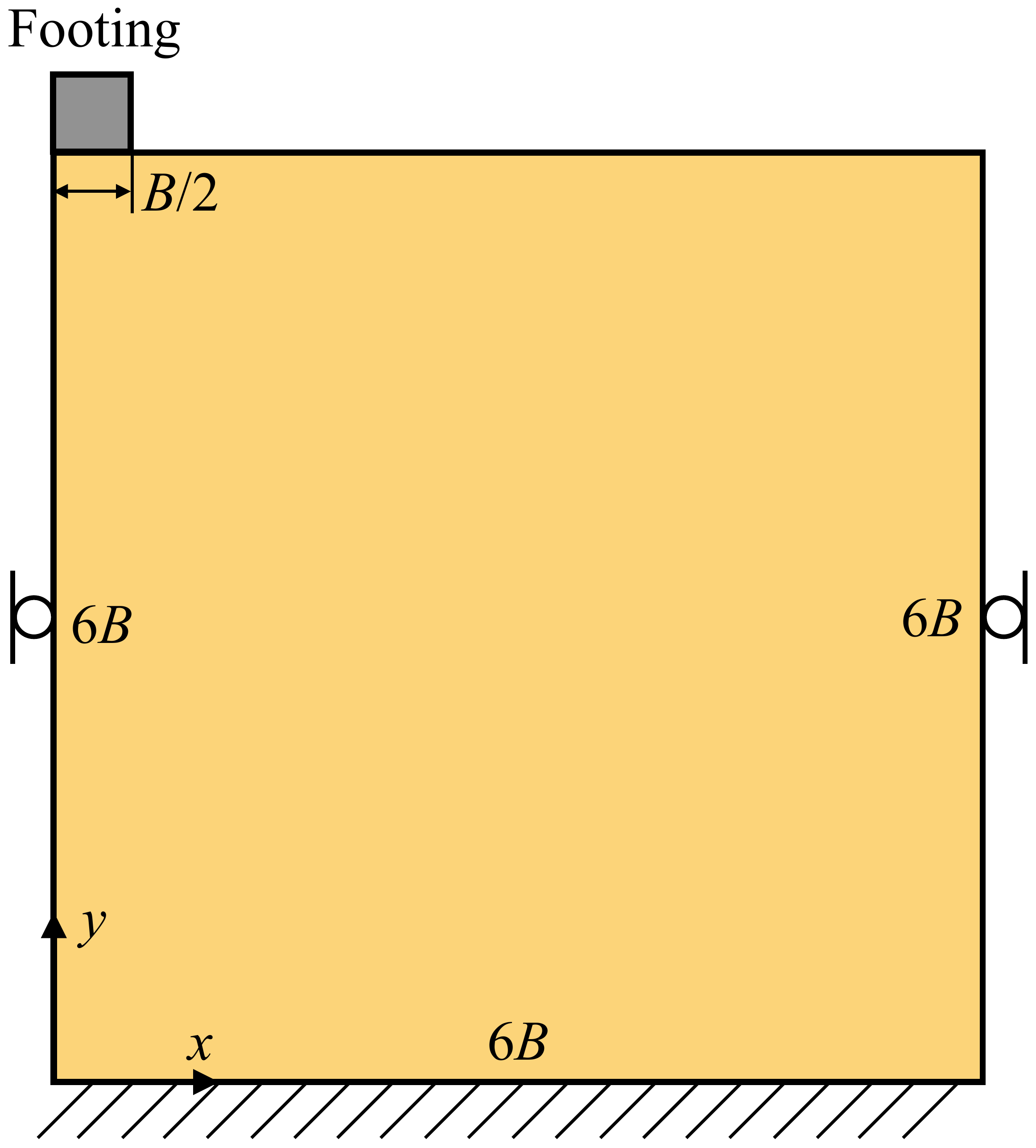}
  \caption{Strip footing: problem geometry and boundary conditions. $B$ refers to the width of the footing.}
  \label{fig:strip_footing_setup}
\end{figure}

In our second example, we investigate the performance of the proposed approach when a nearly incompressible material is in contact with a rigid body---a common scenario in many engineering applications.
To this end, we simulate the problem of a strip footing on an incompressible elastoplastic solid, for which an analytical solution for the bearing capacity---the Prandtl solution---is available.
Figure~\ref{fig:strip_footing_setup} depicts the specific geometry and boundary conditions simulated herein.
As shown in the figure, only the right half of the problem is modeled taking advantage of symmetry. 
Note that the footing is treated explicitly as a rigid body.

We treat the contact between the footing and the ground with a barrier method~\cite{li2020incremental,zhao2022barrier,li2022bfemp,jiang2022hybrid}, which guarantees non-interpenetration between the two objects.
The particular barrier method implemented in this example is based on the formulation tailored to material points in contact with a discrete object~\cite{jiang2022hybrid}.
The friction coefficient in the contact model is set to be sufficiently large to prevent a slip between the footing and the ground.

The elastoplastic behavior of the ground is described by a combination of Hencky elasticity and J2 plasticity.
The elasticity parameters assigned are Young's modulus of $E = 1000$ kPa and Poisson's ratio of $\nu = 0.49$.
The yield strength of J2 plasticity is set such that the (undrained) shear strength of the ground is $0.1$ kPa under plane strain.
These parameters are adopted from Bisht \etal~\cite{bisht2021simulating}.
Note that the high ratio between Young's modulus and the shear strength allows the ground to be in the small deformation range, such that the bearing capacity can be estimated by Prandtl's analytical solution.

Similar to the previous example, we investigate the performance of the proposed approach under three different levels of discretization.
They are:
(i) $h=B/20$ (230,400 material points), 
(ii) $h=B/40$ (921,600 material points), 
and (iii) $h=B/80$ (3,685,400 material points).
The time increment is calculated as $\Delta t = 0.4 (h/c)$, which gives $\Delta t = 1.529 \times 10^{-4}$ s, $\Delta t = 7.644 \times 10^{-5}$ s, and $\Delta t = 3.822 \times 10^{-5}$ s, respectively, for the three levels of grid sizes.
To emulate a quasi-static condition, we apply the damping method used in the previous example (with a damping coefficient of 0.002), as well as pushing the footing slowly with the same penetration rate used in Bisht \etal~\cite{bisht2021simulating}.

\begin{figure}[htbp]
  \centering
  \subfloat[GIMP]{\includegraphics[width=0.45\textwidth]{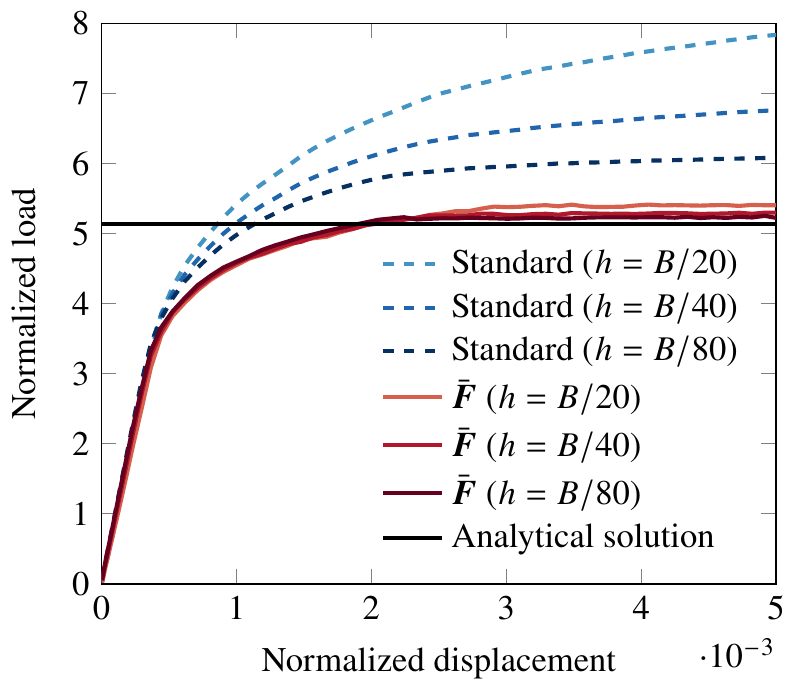}}\hspace{1em}
  \subfloat[B-splines]{\includegraphics[width=0.45\textwidth]{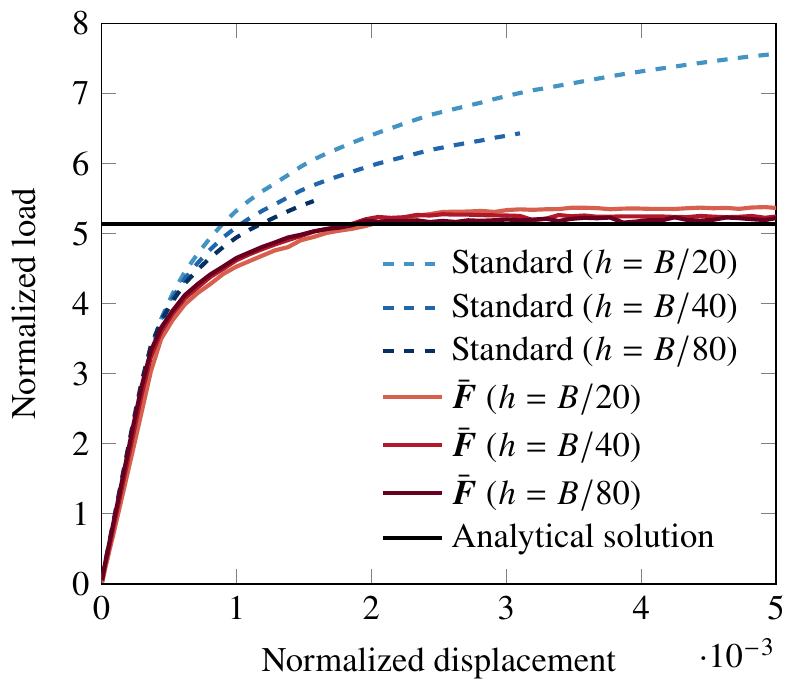}}
  \caption{Strip footing: normalized load--displacement curves from the standard and $\bar{\bm{F}}$ MPM solutions, obtained with GIMP and B-splines basis functions.}
  \label{fig:strip_footing_load_displacement}
\end{figure}

Figure~\ref{fig:strip_footing_load_displacement} presents the normalized load--displacement curves produced from the standard and $\bar{\bm{F}}$ MPM, along with the analytical solution, 5.14. (The load is normalized by the shear strength, and the displacement is normalized by the footing width.)
As is well known, the standard MPM significantly overestimates the bearing capacity due to volumetric locking. 
Also, when B-splines are used, the numerical solutions even blow up, as can be seen from the incomplete load--displacement curves.
By contrast, the $\bar{\bm{F}}$ MPM provides numerical solutions close to the analytical solution.

\begin{figure}[htbp]
  \centering
  \includegraphics[width=0.5\textwidth]{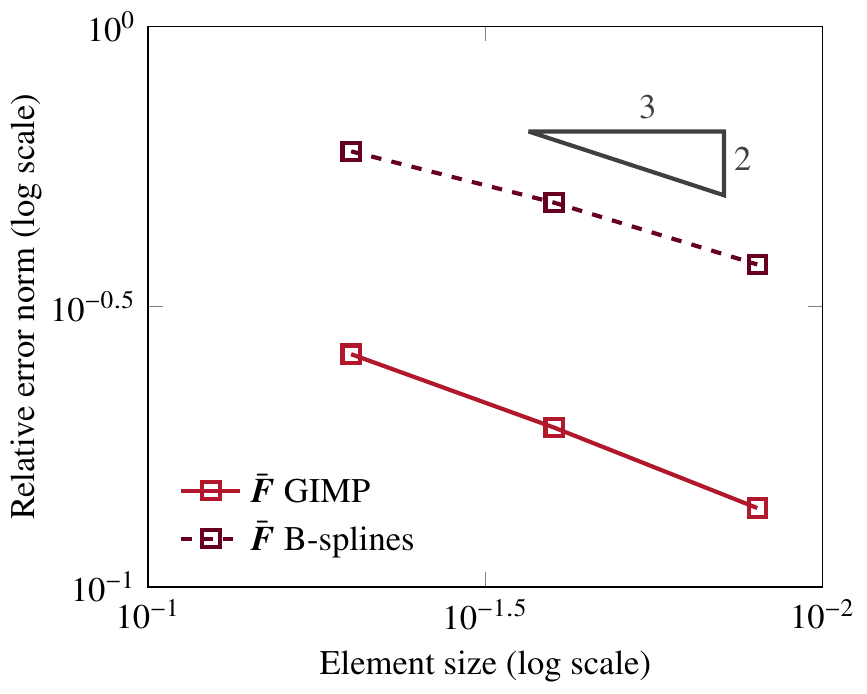}
  \caption{\revised{Strip footing: convergence of the contact pressure solutions with grid refinement.}}
  \label{fig:strip_footing_convergence}
\end{figure}

In Fig.~\ref{fig:strip_footing_convergence} we show the convergence profile of the contact pressure solutions obtained with the $\bar{\bm{F}}$ MPM.
As can be seen, the $\bar{\bm{F}}$ MPM solution converges upon grid refinement, with a rate \revised{of approximately 2/3}.
It is worthwhile to note that the convergence rates of MPM solutions are unavoidably lower than the optimal rates of finite element solutions---even for simple problems with manufactured solutions---due to several reasons including particle--grid transfer operations and sub-optimal quadrature~\revised{\cite{wallstedt2011weighted,kamojjala2015verification,sulsky2016improving}}.  
\revised{Notably, the specific convergence rate in this problem is similar to the convergence rates of MPM solutions to contact problems in the literature (\eg~\cite{deVaucorbeil2021modelling,zhang2021truncated}).
Also, it is observed that the $\bar{\bm{F}}$ MPM solutions obtained with GIMP's basis functions show consistently lower errors than those with B-splines. 
As discussed in Steffen \etal~\cite{steffen2008examination} where the same trend was observed in a numerical example, the main reason would be that the particles remain aligned to the axis in GIMP whereas they do not when B-splines are used.
Overall, the observed convergence rate appears satisfactory.}

\begin{figure}[htbp]
  \centering
  \subfloat[GIMP]{\includegraphics[width=1.0\textwidth]{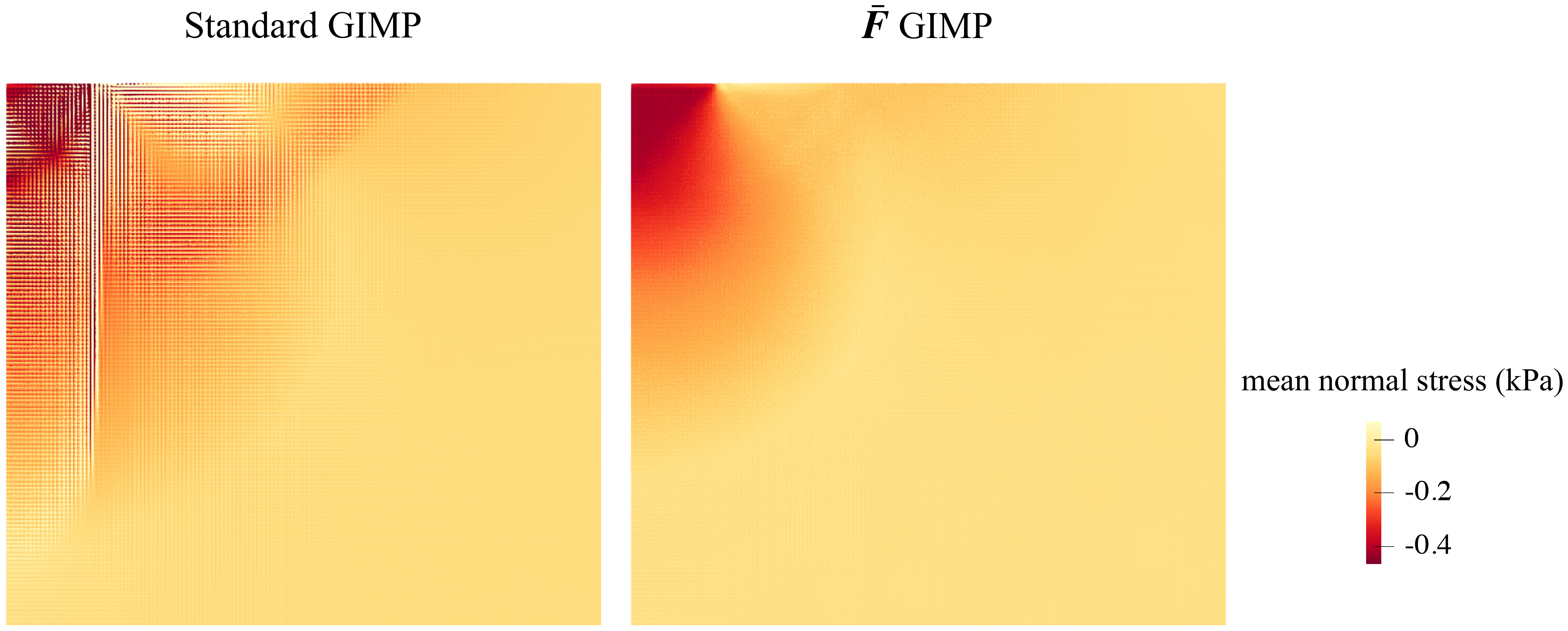}}
  \\
  \subfloat[B-splines]{\includegraphics[width=1.0\textwidth]{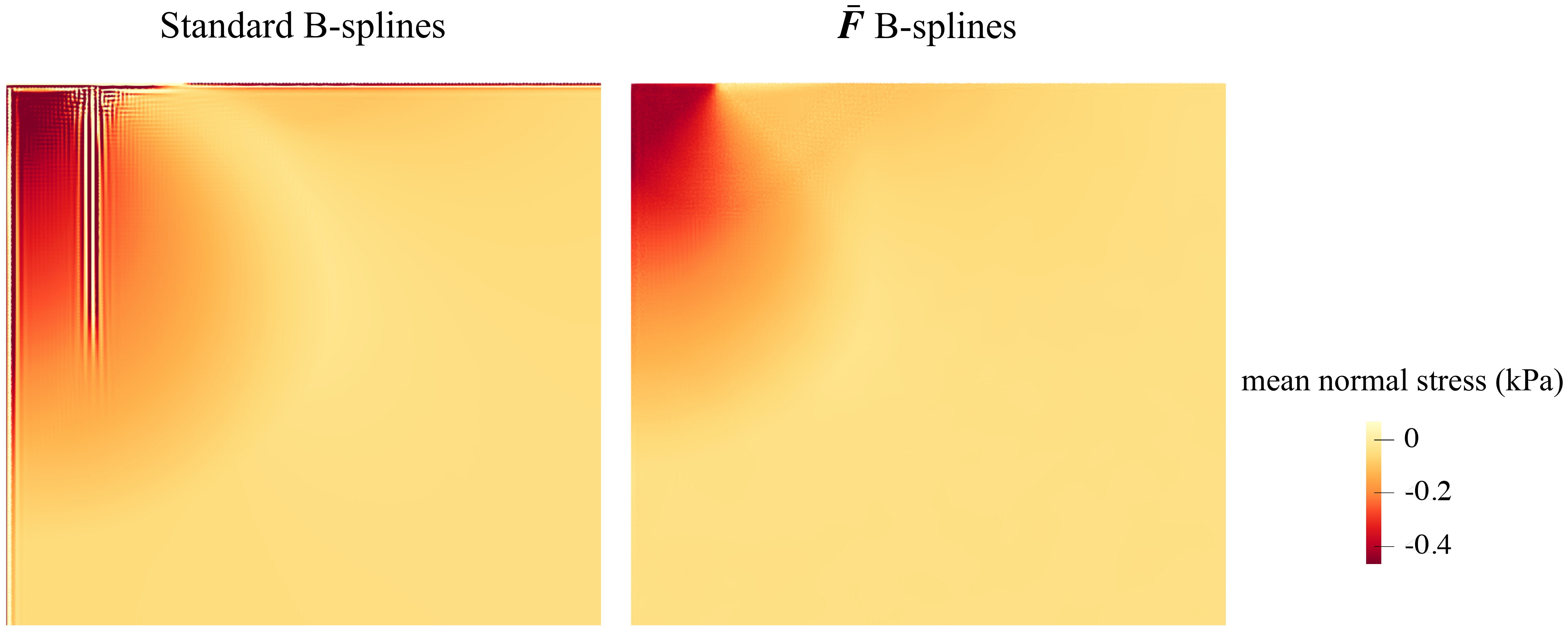}}
  \caption{Strip footing: mean normal stress fields in the standard and $\bar{\bm{F}}$ MPM solutions, obtained with GIMP and B-splines basis functions.}
  \label{fig:strip_footing_mean_normal_stress}
\end{figure}

Figure~\ref{fig:strip_footing_mean_normal_stress} shows the mean normal stress fields in the standard and $\bar{\bm{F}}$ MPM solutions when $h=B/40$.
As in Cook's membrane example, the standard MPM solutions show non-physical oscillations in the stress fields, which are particularly severe when GIMP is used.
The solutions obtained by the proposed $\bar{\bm{F}}$ MPM, however, are free of such oscillations.
Considering this result together with the bearing capacity results above, it can be concluded that the proposed method successfully alleviates volumetric locking in the MPM involving contact.

\subsection{Dam break}
\begin{figure}[htbp]
  \centering
  \includegraphics[width=0.65\textwidth]{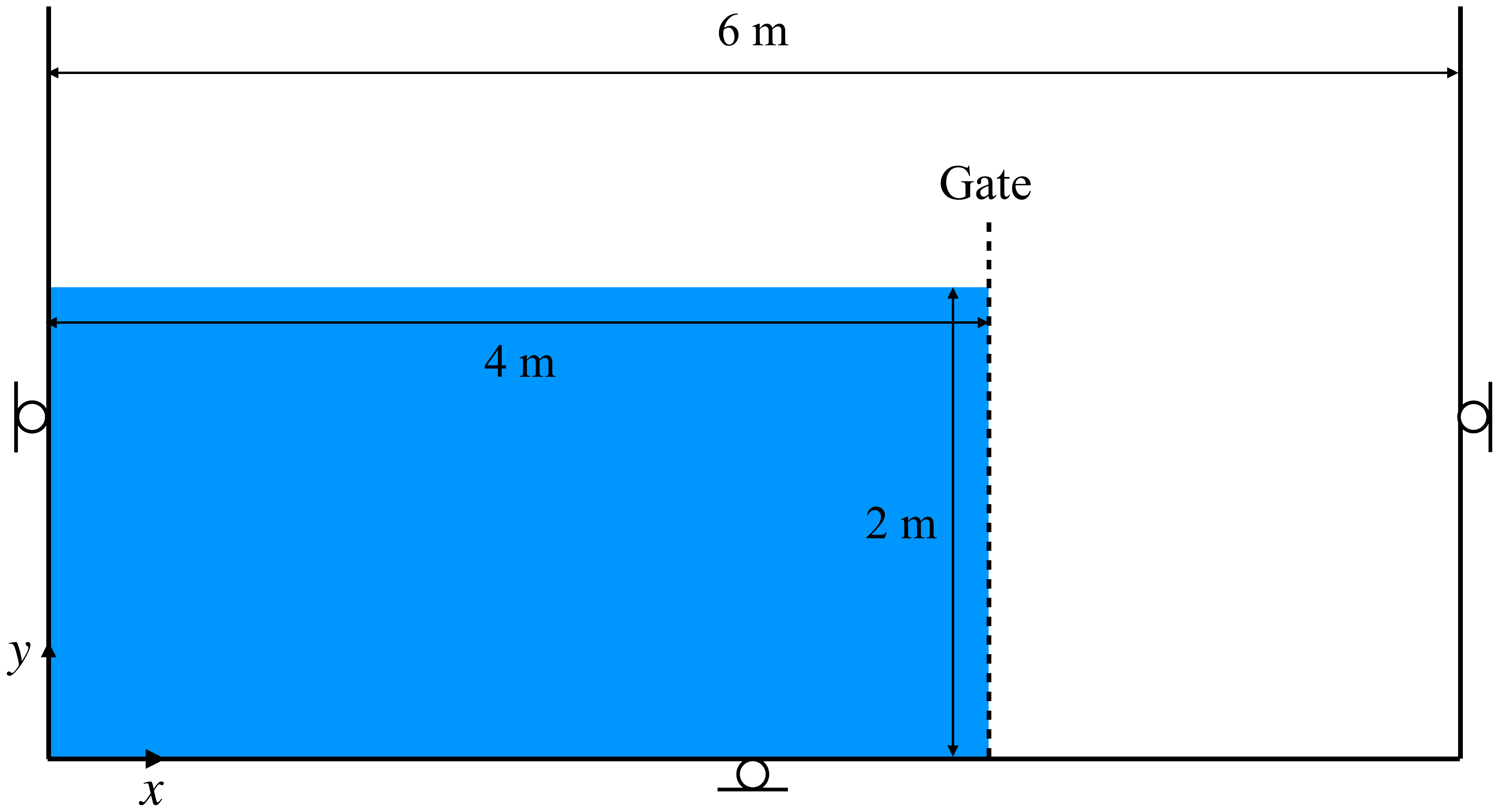}
  \caption{Dam break: problem geometry and boundary conditions.}
  \label{fig:dam_break_setup}
\end{figure}

Our third example is the dam break problem in Mast \etal~\cite{mast2012mitigating}, where a three-field mixed formulation is used for mitigating locking in the MPM.
Figure~\ref{fig:dam_break_setup} depicts the geometry and boundary conditions of the problem.
As shown, it considers a 4-m-long and 2-m-high water reservoir initially constrained by a gate and starts to flow after the gate is removed. 
It is noted that there is another gate on the right boundary, from which the water will bounce back.

Following Mast \etal~\cite{mast2012mitigating}, the water is modeled as a nearly incompressible Newtonian fluid, with a bulk modulus of $K = 2.0$ GPa, a dynamic viscosity of $\mu = 0.001$ Pa$\cdot$s, and a density of $\rho$ = 0.9975 t/m$^3$.
For MPM discretization, we introduce a background grid comprised of 0.25-m long square elements and initialize each element with 25 material points. 
This results in 3,200 material points in total.
We then apply gravity loading to the water reservoir until it reaches the hydrostatic state. 
Subsequently, we remove the gate so that the water can flow freely.
We simulate the problem until $t=2.0$ s, setting the time increment as $\Delta t = 10^{-5}$ s. 
    
\begin{figure}
  \centering
  \subfloat[GIMP]{\includegraphics[width=1.0\textwidth]{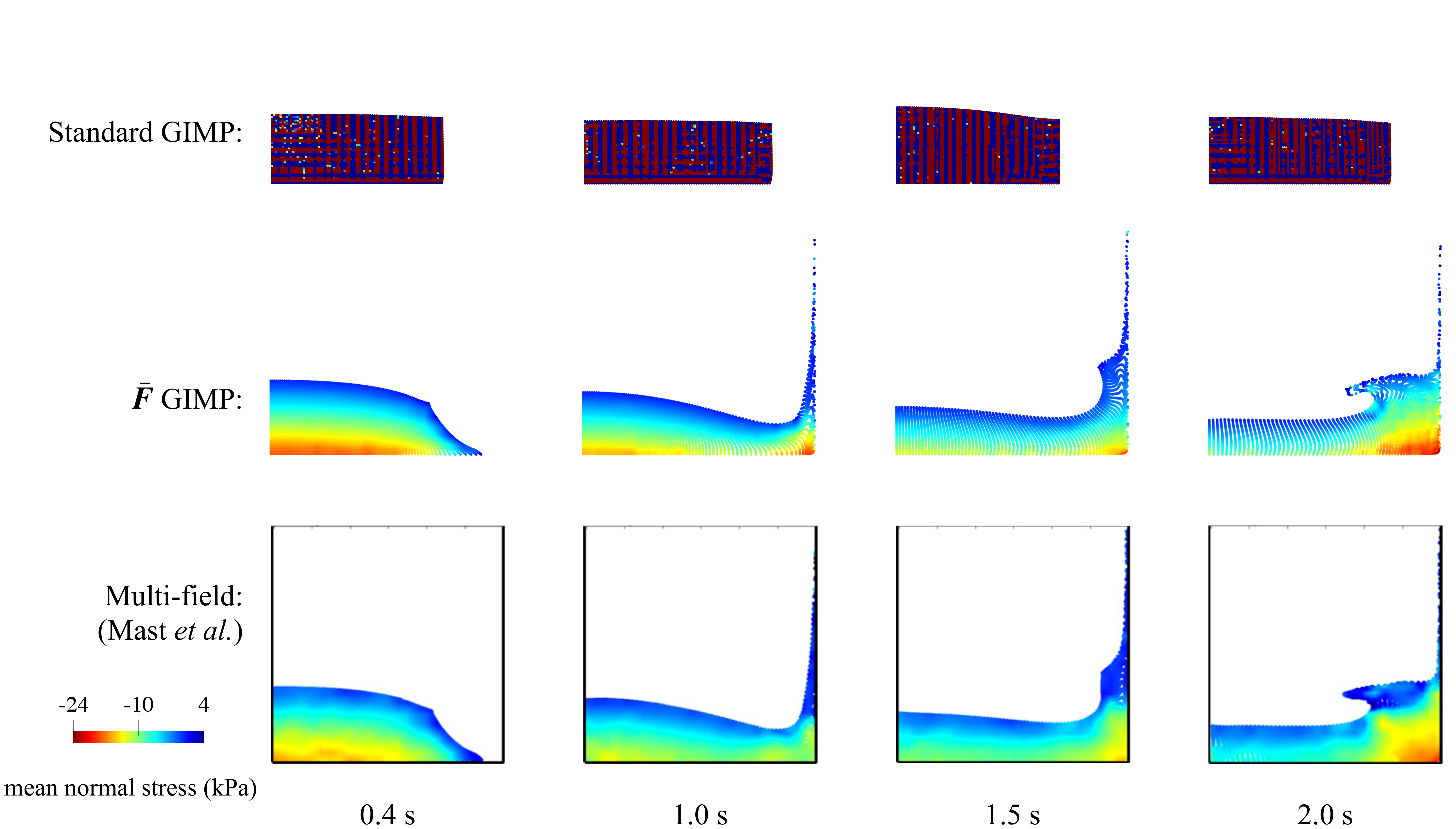}}\\
  \subfloat[B-splines]{\includegraphics[width=1.0\textwidth]{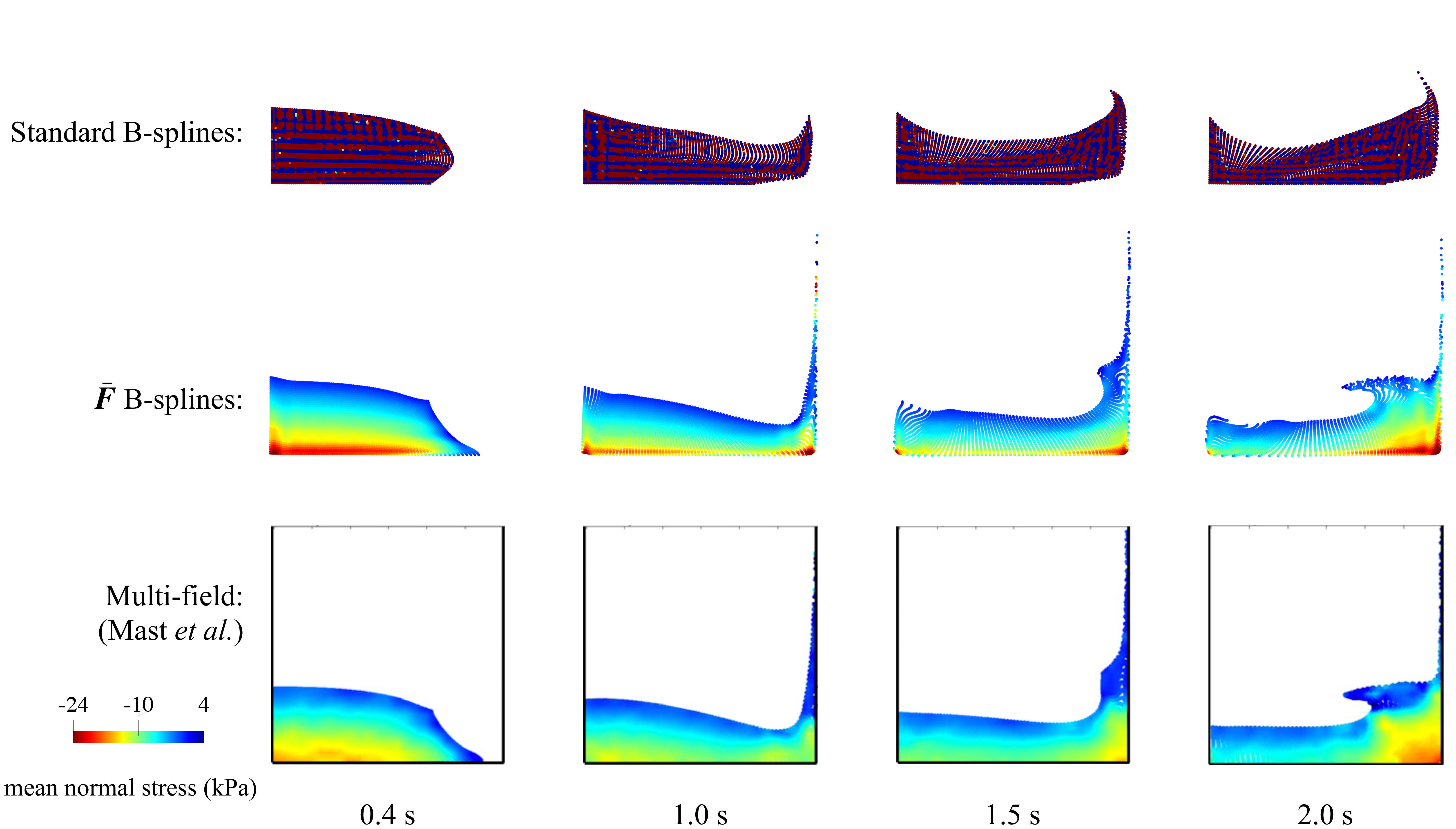}}
  \caption{Dam break: flow snapshots from the standard and $\bar{\bm{F}}$ MPM simulations, along with those from Mast \etal~\cite{mast2012mitigating} produced from a three-field mixed MPM formulation. Material points are colored by the water pressure.}
  \label{fig:dam_break_verification}
\end{figure}

Figure~\ref{fig:dam_break_verification} shows the flow snapshots simulated by the standard and $\bar{\bm{F}}$ MPM formulations (with both GIMP and B-splines), in comparison with the reference solutions from Mast \etal~\cite{mast2012mitigating}.
One can easily see that the standard MPM is subjected to severe volumetric locking.
Not only does the pressure field show non-physical oscillations, but the water also is unrealistically stiff. 
The locking is less severe in the B-splines MPM than GIMP, but it is still unacceptable.
When the $\bar{\bm{F}}$ MPM formulation is employed, however, the numerical solutions are free of the locking problem, irrespective of the basis functions.
It can also be seen that the numerical solutions produced by the $\bar{\bm{F}}$ MPM formulation are very similar to those by the three-field mixed MPM formulation and linear shape functions~\cite{mast2012mitigating}. 
\revisedsecond{Note that while the B-splines MPM solutions show some difference around the boundaries, it is known that the B-splines MPM has issues near boundaries, which has been the subject of other works (\eg~\cite{schulz2019consistent,nakamura2023taylor}). }
Importantly, it is reminded that the proposed method involves significantly less implementation effort and computational cost than the three-field mixed formulation, while the results are nearly the same. 
It is also worthwhile to note that the applicability of the three-field mixed formulation is tied to that of the Hu--Washizu variational principle.

\begin{figure}
  \centering
  \includegraphics[width=0.5\textwidth]{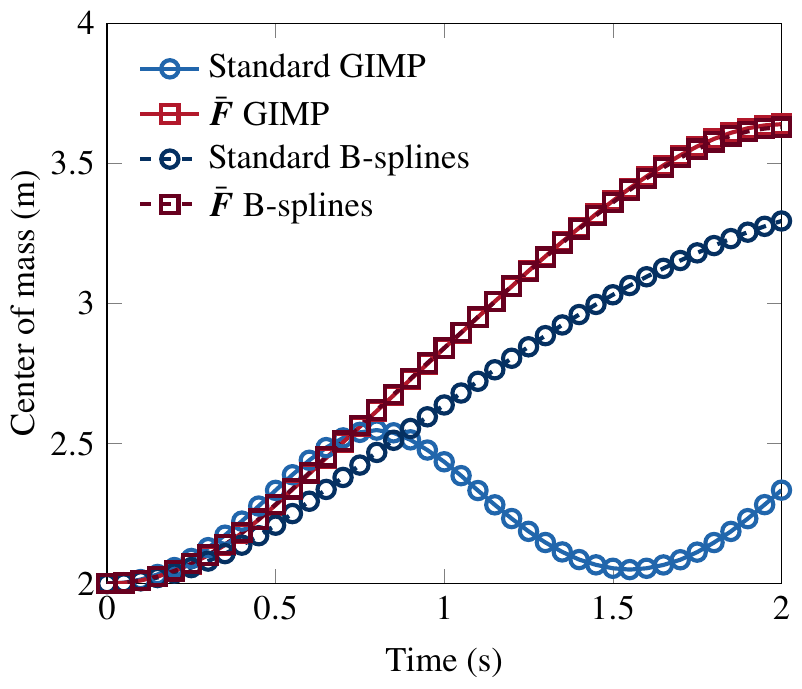}
  \caption{\revised{Dam break: time evolution of the center-of-mass in the flow direction in the standard and $\bar{\tensor{F}}$ MPM solutions, obtained with GIMP and B-splines basis functions.}}
  \label{fig:dam_break_center_of_mass_position}
\end{figure}

\revised{Further, for a quantitative comparison, in Fig.~\ref{fig:dam_break_center_of_mass_position} we plot the time evolutions of the center-of-mass in the flow direction calculated from the standard and $\bar{\tensor{F}}$ MPM solutions, similar to an analysis performed in Baumgarten and Kamrin~\cite{baumgarten2023analysis}.
It can be seen that the center-of-mass evolves nearly identically in the $\bar{\tensor{F}}$ MPM formulations obtained with GIMP and B-splines basis function, while the standard MPM solutions exhibit significant differences.  
This comparison confirms the proposed method performs well and consistently irrespective of the basis functions.
}

\subsection{3D landslide}
\begin{figure}
  \centering
  \includegraphics[width=0.5\textwidth]{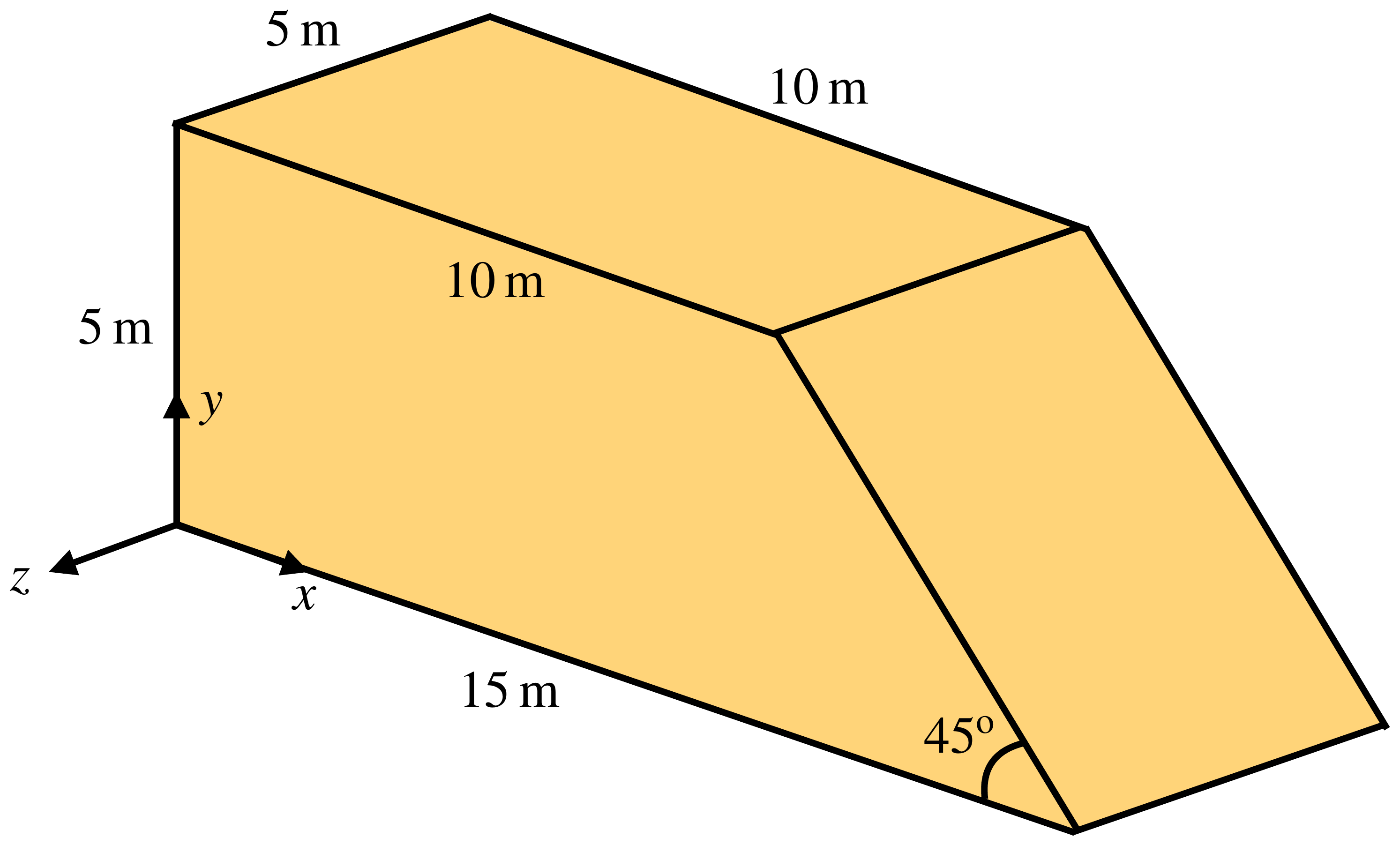}
  \caption{3D landslide: problem geometry.}
  \label{fig:slope_setup}
\end{figure}

As our last example, we investigate the performance of the proposed approach for 3D large deformation in an elastoplastic material. 
To this end, we simulate a 3D landslide process where a brittle clay slope fails in an undrained manner.
The slope geometry is depicted in Fig.~\ref{fig:slope_setup}.
The bottom boundary of the slope is fully fixed, while the three lateral boundaries are supported by rollers.
The elastoplastic behavior of the undrained clay is modeled through a combination of Hencky elasticity and J2 plasticity with a softening law. 
The specific softening law adopted is $\kappa = \kappa_r + (\kappa_p - \kappa_r) e^{-\eta \varepsilon_q^{\rm p}}$, where $\kappa$ is the yield strength, $\kappa_r$ and $\kappa_p$ are the residual and peak strengths, respectively, $\eta$ is the softening parameter, and $\varepsilon_q^{\rm p}$ is the cumulative equivalent plastic strain.
The material parameters are assigned to be similar to the undrained sensitive clay modeled in Bui and Nguyen~\cite{bui2021smoothed}.
They are: Young's modulus $E = 25$ MPa, Poisson's ratio $\nu = 0.499$, peak strength $\kappa_p = 40.82$ kPa, residual strength $\kappa_r = 2.45$ kPa, softening coefficient $\eta = 5$, and the density $\rho = 2.15$ t/m$^3$.

For MPM discretization, we introduce a 3D background grid comprised of mono-sized cubic elements whose length is 0.2 m. 
For each element in the slope domain, we assign 8 material points, which results in a total of 311,250 material points.
We initialize the stress field in the material points through a gravity loading stage.  
Then, to trigger the slope failure, we decrease the peak strength with a reduction factor of 1.65, as done in Bui and Nguyen~\cite{bui2021smoothed}.
We simulate the problem until $5.5$ s with a time increment of $\Delta t = 5 \times 10^{-5}$ s.

\begin{figure}[htbp]
  \centering
  \subfloat[$t = 1.5$ s]{\includegraphics[width=0.85\textwidth]{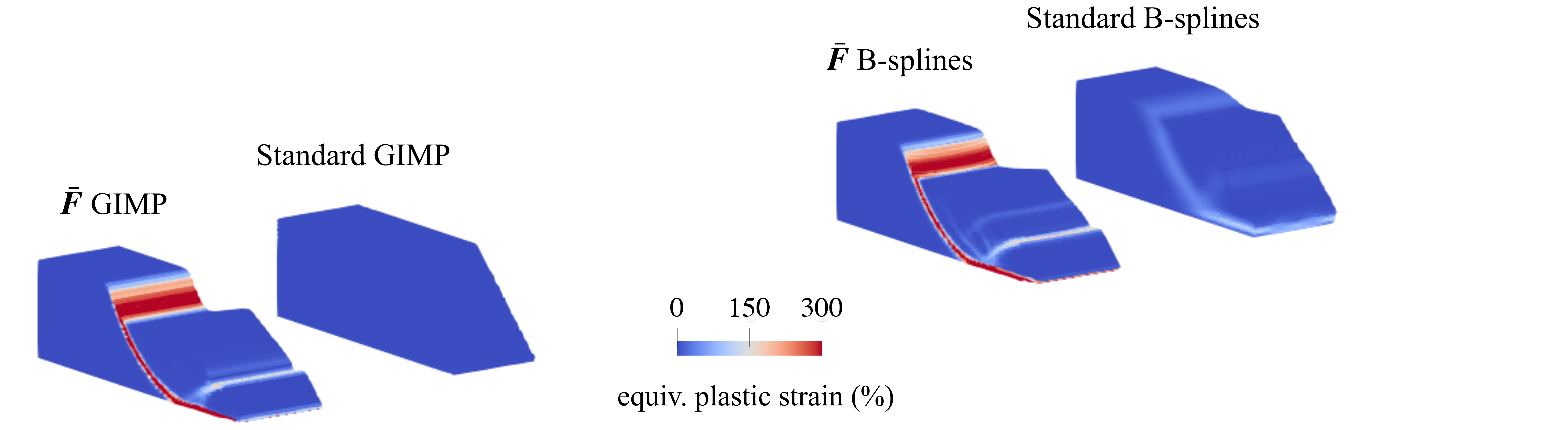}}\\
  \subfloat[$t = 2.5$ s]{\includegraphics[width=0.85\textwidth]{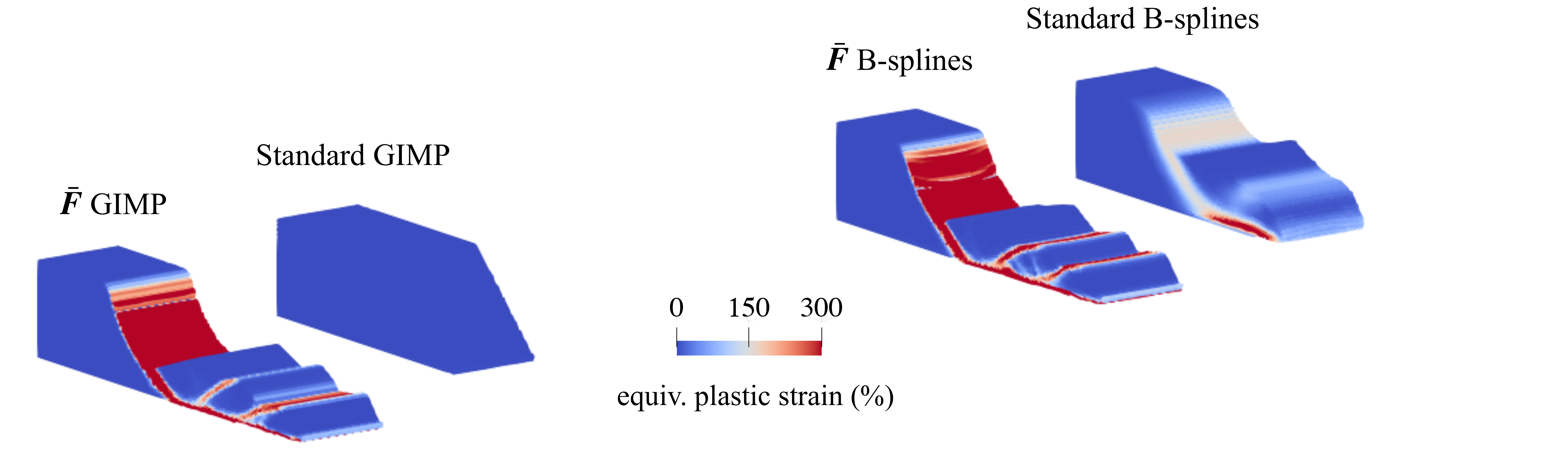}}\\
  \subfloat[$t = 3.5$ s]{\includegraphics[width=0.85\textwidth]{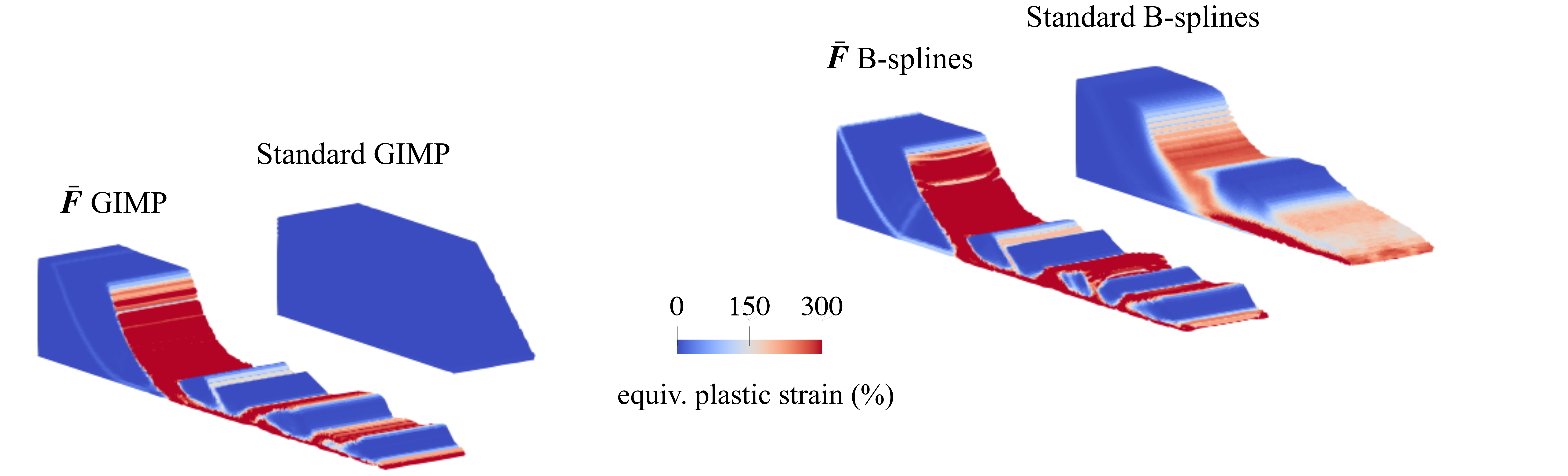}}\\
  \subfloat[$t = 5.5$ s]{\includegraphics[width=0.85\textwidth]{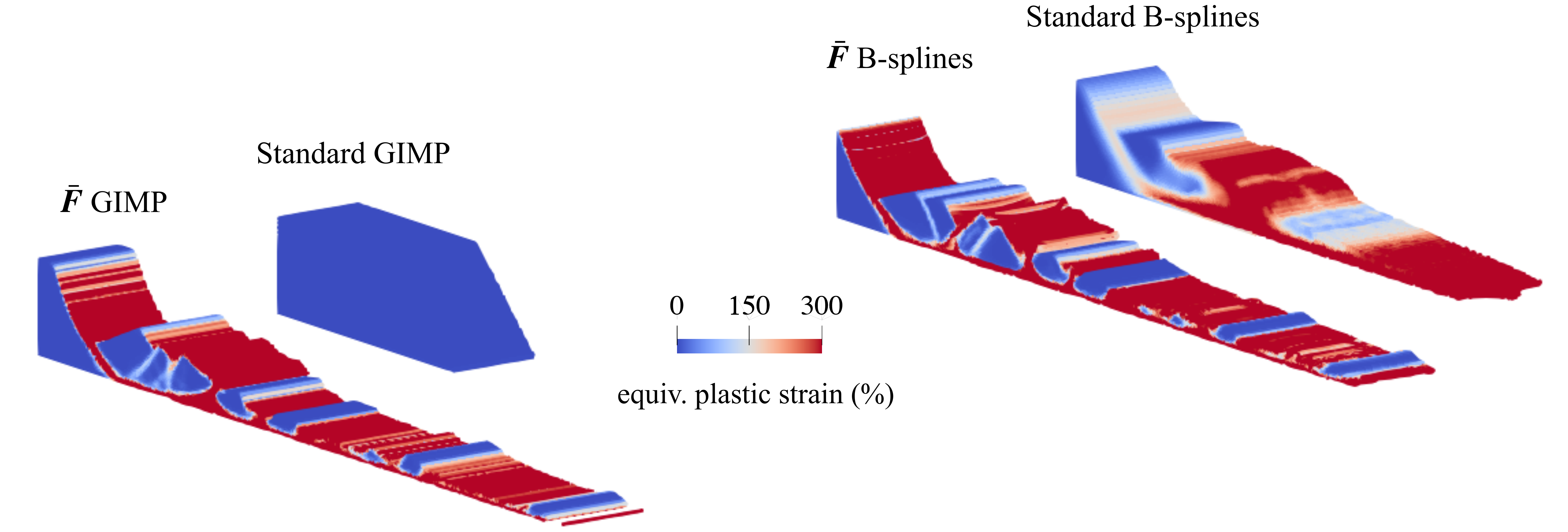}}
  \caption{3D landslides: snapshots from the standard and $\bar{\bm{F}}$ MPM simulations. Material points are colored by the equivalent plastic strain.}
  \label{fig:slope_strain}
\end{figure}

\begin{figure}[htbp]
  \centering
  \subfloat[$t = 1.5$ s]{\includegraphics[width=0.85\textwidth]{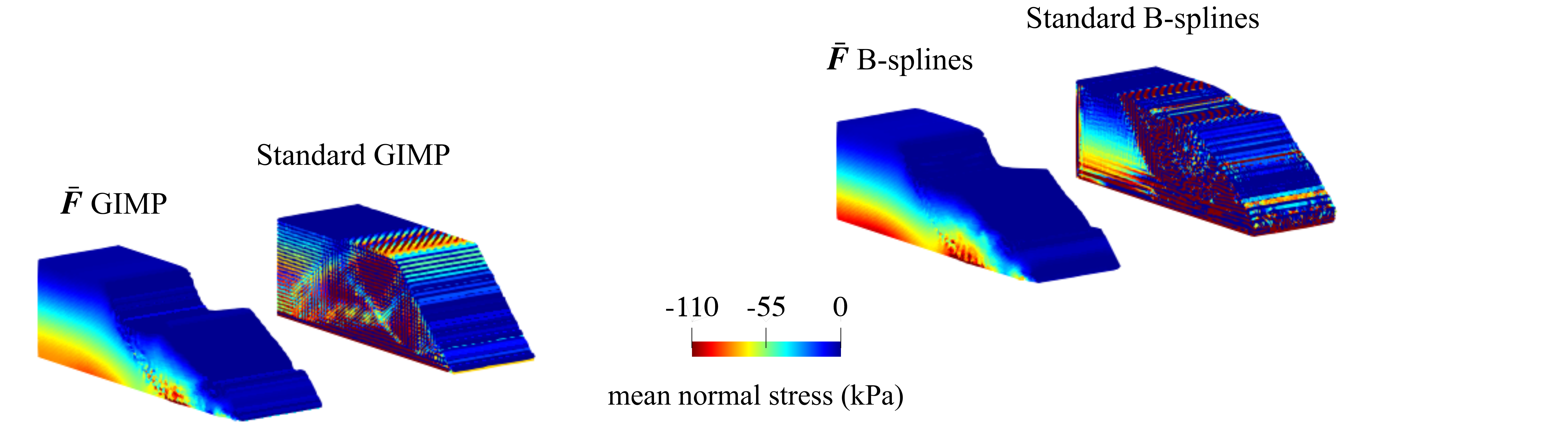}}\\
  \subfloat[$t = 2.5$ s]{\includegraphics[width=0.85\textwidth]{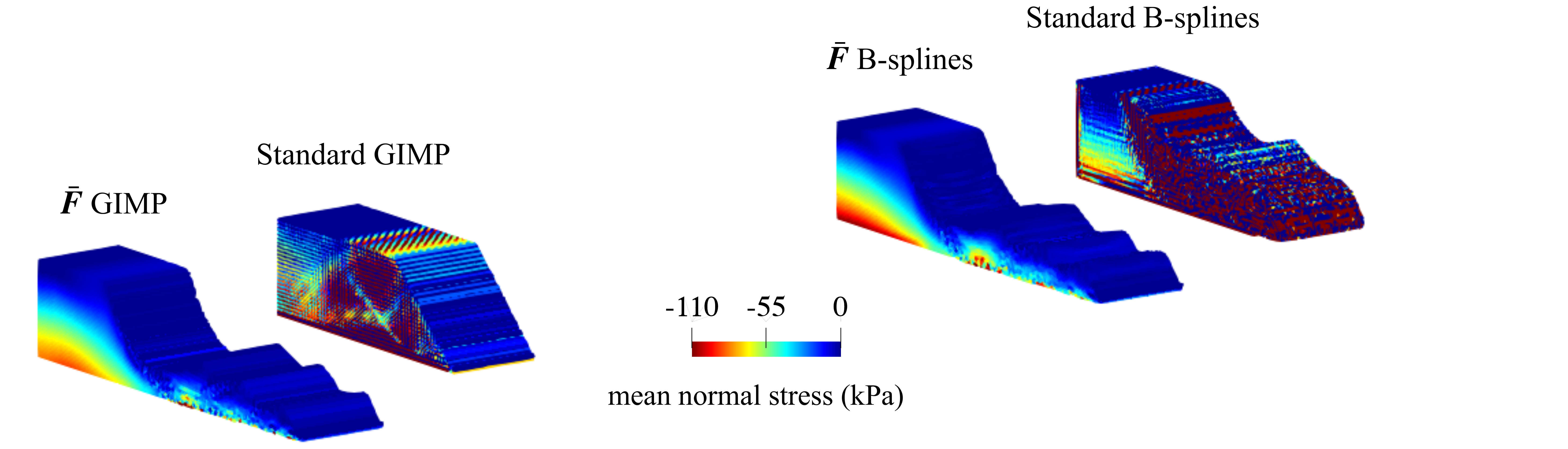}}\\
  \subfloat[$t = 3.5$ s]{\includegraphics[width=0.85\textwidth]{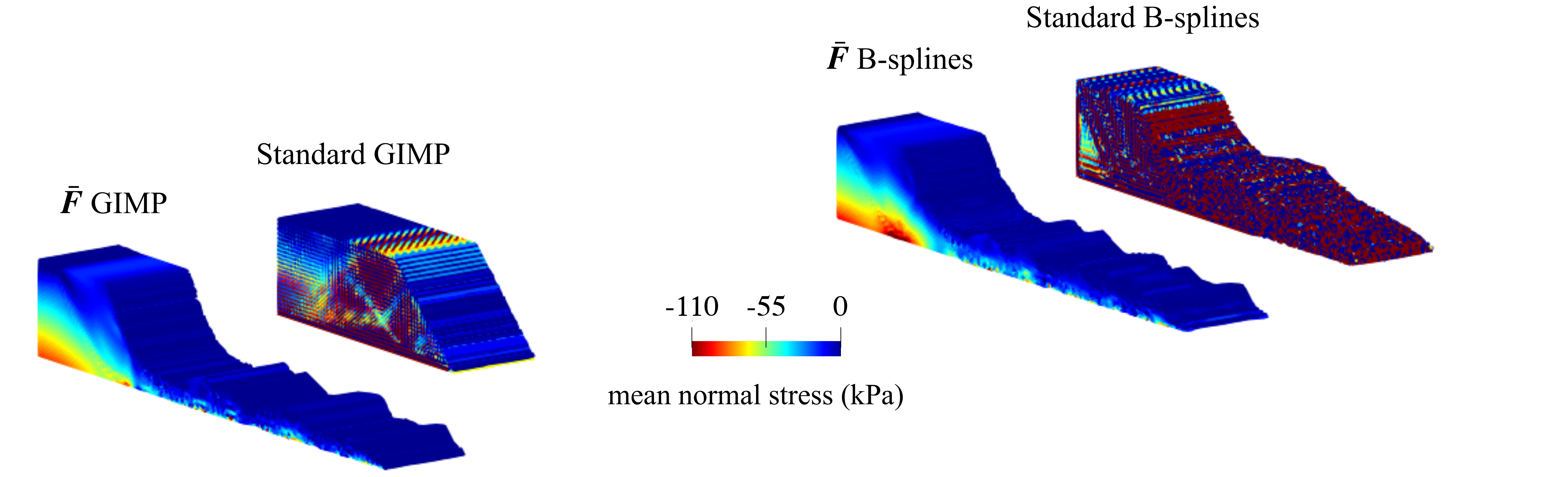}}\\
  \subfloat[$t = 5.5$ s]{\includegraphics[width=0.85\textwidth]{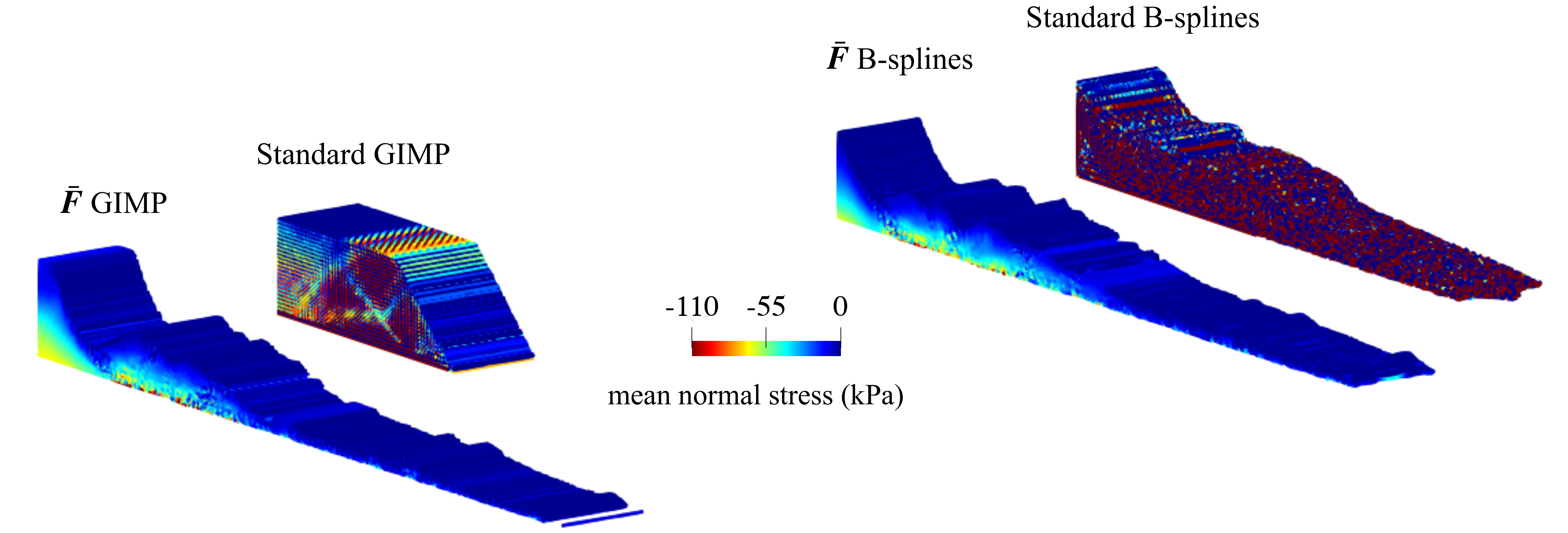}}
  \caption{3D landslides: snapshots from the standard and $\bar{\bm{F}}$ MPM simulations. Material points are colored by the mean normal stress.}
  \label{fig:slope_stress}
\end{figure}

Figures~\ref{fig:slope_strain} and~\ref{fig:slope_stress} present snapshots of the landslide simulated by the standard and $\bar{\bm{F}}$ MPM, showing the equivalent plastic strain and mean normal stress fields, respectively.
As in the previous example, the standard MPM is subjected to severe volumetric locking, manifesting non-physical stress oscillations as well as overly stiff behavior.
The slope does not even fail when GIMP is used, while it shows diffusive plastic strains when B-splines are used. 
However, the $\bar{\bm{F}}$ MPM results show retrogressive failure---a signature failure pattern in sensitive clay slopes---with both GIMP's basis functions and B-splines.
Although there is no reference solution for this problem, it can be seen that the simulation results produced with the two types of basis functions are quite similar.
Thus the $\bar{\bm{F}}$ MPM solutions are believed to be reasonable.

\begin{figure}[htbp]
  \centering
  \subfloat[GIMP]{\includegraphics[width=1.0\textwidth]{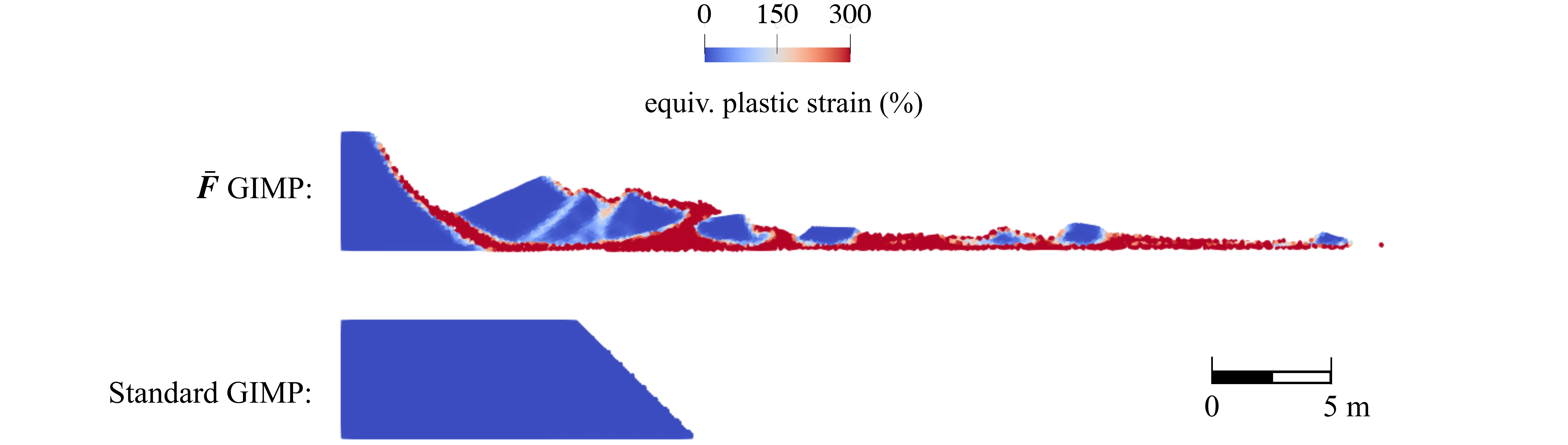}}\\
  \subfloat[B-splines]{\includegraphics[width=1.0\textwidth]{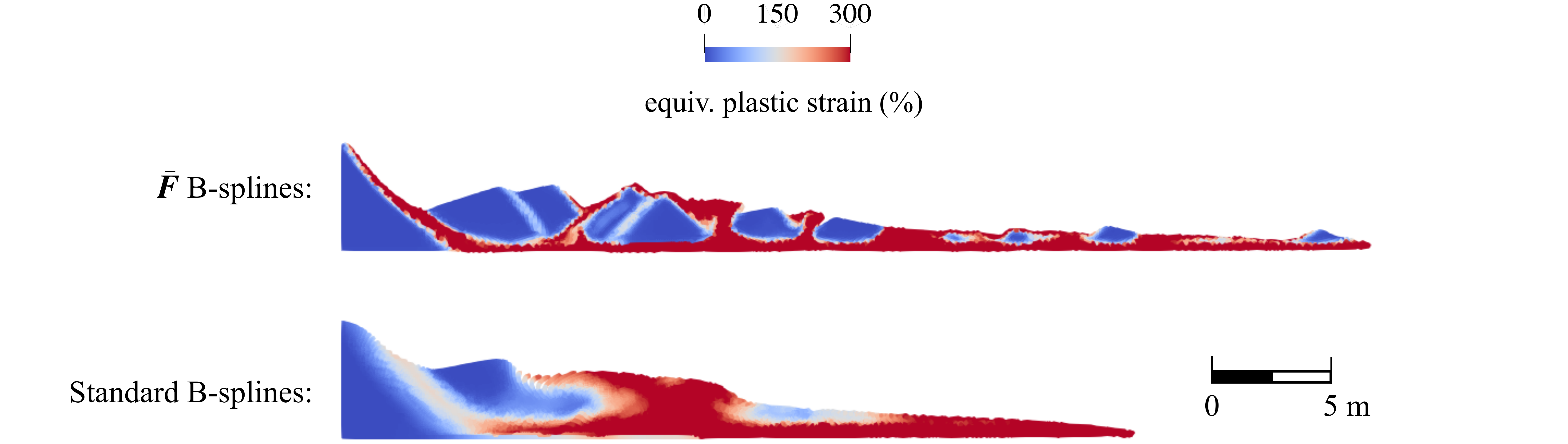}}
  \caption{\revised{3D landslides: comparison of run-out distances in the standard and $\bar{\bm{F}}$ MPM simulations.}}
  \label{fig:slope_runout}
\end{figure}

Figure~\ref{fig:slope_runout} demonstrates how the run-out distance, which is of primary interest in slope analysis, is different in the standard and $\bar{\bm{F}}$ MPM solutions.
\revisedsecond{The run-out distances are: 0.0055 m (standard GIMP), 19.07 m (standard B-splines), 29.96 m ($\bar{\bm{F}}$ GIMP), and 29.29 m ($\bar{\bm{F}}$ B-splines).}
Without a proper locking-mitigation approach, the standard MPM significantly underestimates the run-out distance, if not being unable to simulate the failure process at all. 
This difference highlights why it is of critical importance to mitigate volumetric locking from the practical viewpoint.


\section{Closure}
\label{section:closure}
  
In this paper, we have proposed a simple and efficient approach for circumventing volumetric locking in the MPM, which can be generally applied to a family of standard explicit MPM formulations regardless of basis functions and material types.
The key idea of the proposed approach is to evaluate the assumed deformation gradient ($\bar{\bm{F}}$) with a volume-averaging operation that resembles the standard particle--grid transfer scheme in the MPM.
The approach can be implemented in a much simpler way than the existing approaches for mitigating locking in the MPM, and it is independent of other parts such as the basis function and the constitutive behavior.
The results of the numerical examples have verified and demonstrated that the proposed approach performs well for various types of nearly incompressible problems arising in solid and fluid mechanics.
Taken together, it is believed that the proposed approach is highly attractive for mitigating volumetric locking in the explicit MPM.

\section*{Acknowledgments}

This work was supported by the National Research Foundation of Korea (NRF) grant funded by the Korean government (MSIT) (Nos. 2022R1F1A1065418 and RS-2023-00209799). 
The authors also wish to thank Dr. Vibhav Bisht for sharing his numerical solution to Cook's membrane problem obtained with the nonlinear $\bar{\bm{B}}$ MPM.

\section*{Data Availability Statement} 
\label{sec:data-availability} 

The data that support the findings of this study are available from the corresponding author upon reasonable request.

\bibliography{references}

\begin{thebibliography}{10}
\expandafter\ifx\csname url\endcsname\relax
  \def\url#1{\texttt{#1}}\fi
\expandafter\ifx\csname urlprefix\endcsname\relax\def\urlprefix{URL }\fi
\expandafter\ifx\csname href\endcsname\relax
  \def\href#1#2{#2} \def\path#1{#1}\fi

\bibitem{sulsky1994particle}
D.~Sulsky, Z.~Chen, H.~L. Schreyer, A particle method for history-dependent
  materials, Computer Methods in Applied Mechanics and Engineering 118~(1-2)
  (1994) 179--196.

\bibitem{brackbill1986flip}
J.~U. Brackbill, H.~M. Ruppel, {FLIP: A method for adaptively zoned,
  particle-in-cell calculations of fluid flows in two dimensions}, Journal of
  Computational Physics 65~(2) (1986) 314--343.

\bibitem{gaume2018dynamic}
J.~Gaume, T.~Gast, J.~Teran, A.~Van~Herwijnen, C.~Jiang, Dynamic anticrack
  propagation in snow, Nature Communications 9~(1) (2018) 1--10.

\bibitem{fern2019material}
J.~Fern, A.~Rohe, K.~Soga, E.~Alonso, The Material Point Method for
  Geotechnical Engineering: A Practical Guide, CRC Press, 2019.

\bibitem{zhao2020stabilized}
Y.~Zhao, J.~Choo, Stabilized material point methods for coupled large
  deformation and fluid flow in porous materials, Computer Methods in Applied
  Mechanics and Engineering 362 (2020) 112742.

\bibitem{li2021three}
X.~Li, B.~Sovilla, C.~Jiang, J.~Gaume, Three-dimensional and real-scale
  modeling of flow regimes in dense snow avalanches, Landslides (2021) 1--14.

\bibitem{love2006energy}
E.~Love, D.~L. Sulsky, An energy-consistent material-point method for dynamic
  finite deformation plasticity, International Journal for Numerical Methods in
  Engineering 65~(10) (2006) 1608--1638.

\bibitem{mast2012mitigating}
C.~Mast, P.~Mackenzie-Helnwein, P.~Arduino, G.~R. Miller, W.~Shin, Mitigating
  kinematic locking in the material point method, Journal of Computational
  Physics 231~(16) (2012) 5351--5373.

\bibitem{iaconeta2019stabilized}
I.~Iaconeta, A.~Larese, R.~Rossi, E.~O{\~n}ate, A stabilized mixed implicit
  material point method for non-linear incompressible solid mechanics,
  Computational Mechanics 63~(6) (2019) 1243--1260.

\bibitem{zhang2017incompressible}
F.~Zhang, X.~Zhang, K.~Y. Sze, Y.~Lian, Y.~Liu, Incompressible material point
  method for free surface flow, Journal of Computational Physics 330 (2017)
  92--110.

\bibitem{kularathna2017implicit}
S.~Kularathna, K.~Soga, Implicit formulation of material point method for
  analysis of incompressible materials, Computer Methods in Applied Mechanics
  and Engineering 313 (2017) 673--686.

\bibitem{chorin1968numerical}
A.~J. Chorin, Numerical solution of the {Navier--Stokes} equations, Mathematics
  of Computation 22~(104) (1968) 745--762.

\bibitem{de1996design}
E.~de~Souza~Neto, D.~Peri{\'c}, M.~Dutko, D.~Owen, Design of simple low order
  finite elements for large strain analysis of nearly incompressible solids,
  International Journal of Solids and Structures 33~(20-22) (1996) 3277--3296.

\bibitem{coombs2018overcoming}
W.~M. Coombs, T.~J. Charlton, M.~Cortis, C.~E. Augarde, Overcoming volumetric
  locking in material point methods, Computer Methods in Applied Mechanics and
  Engineering 333 (2018) 1--21.

\bibitem{wang2021efficient}
L.~Wang, W.~M. Coombs, C.~E. Augarde, M.~Cortis, M.~J. Brown, A.~J. Brennan,
  J.~A. Knappett, C.~Davidson, D.~Richards, D.~J. White, et~al., An efficient
  and locking-free material point method for three-dimensional analysis with
  simplex elements, International Journal for Numerical Methods in Engineering
  (2021).

\bibitem{moutsanidis2020treatment}
G.~Moutsanidis, J.~J. Koester, M.~R. Tupek, J.-S. Chen, Y.~Bazilevs, Treatment
  of near-incompressibility in meshfree and immersed-particle methods,
  Computational Particle Mechanics 7~(2) (2020) 309--327.

\bibitem{simo1985variational}
J.~Simo, R.~L. Taylor, K.~Pister, Variational and projection methods for the
  volume constraint in finite deformation elasto-plasticity, Computer Methods
  in Applied Mechanics and Engineering 51~(1-3) (1985) 177--208.

\bibitem{hughes1980generalization}
T.~J. Hughes, Generalization of selective integration procedures to anisotropic
  and nonlinear media, International Journal for Numerical Methods in
  Engineering 15~(9) (1980) 1413--1418.

\bibitem{bisht2021simulating}
V.~Bisht, R.~Salgado, M.~Prezzi, Simulating penetration problems in
  incompressible materials using the material point method, Computers and
  Geotechnics 133 (2021) 103593.

\bibitem{bisht2021material}
V.~Bisht, R.~Salgado, M.~Prezzi, Material point method for cone penetration in
  clays, Journal of Geotechnical and Geoenvironmental Engineering 147~(12)
  (2021) 04021158.

\bibitem{telikicherla2022treatment}
R.~M. Telikicherla, G.~Moutsanidis, Treatment of near-incompressibility and
  volumetric locking in higher order material point methods, Computer Methods
  in Applied Mechanics and Engineering 395 (2022) 114985.

\bibitem{steffen2008analysis}
M.~Steffen, R.~M. Kirby, M.~Berzins, {Analysis and reduction of quadrature
  errors in the material point method (MPM)}, International Journal for
  Numerical Methods in Engineering 76~(6) (2008) 922--948.

\bibitem{gan2018enhancement}
Y.~Gan, Z.~Sun, Z.~Chen, X.~Zhang, Y.~Liu, {Enhancement of the material point
  method using B-spline basis functions}, International Journal for Numerical
  Methods in Engineering 113~(3) (2018) 411--431.

\bibitem{jiang2016material}
C.~Jiang, C.~Schroeder, J.~Teran, A.~Stomakhin, A.~Selle, The material point
  method for simulating continuum materials, in: ACM SIGGRAPH 2016 Courses,
  2016, pp. 1--52.

\bibitem{zhang2016material}
X.~Zhang, Z.~Chen, Y.~Liu, The Material Point Method: A Continuum-Based
  Particle Method for Extreme Loading Cases, Academic Press, 2016.

\bibitem{de2020material}
A.~de~Vaucorbeil, V.~P. Nguyen, S.~Sinaie, J.~Y. Wu, Material point method
  after 25 years: Theory, implementation, and applications, Advances in Applied
  Mechanics 53 (2020) 185--398.

\bibitem{jiang2015affine}
C.~Jiang, C.~Schroeder, A.~Selle, J.~Teran, A.~Stomakhin, The affine
  particle-in-cell method, ACM Transactions on Graphics (TOG) 34~(4) (2015)
  1--10.

\bibitem{fu2017polynomial}
C.~Fu, Q.~Guo, T.~Gast, C.~Jiang, J.~Teran, A polynomial particle-in-cell
  method, ACM Transactions on Graphics (TOG) 36~(6) (2017) 1--12.

\bibitem{hammerquist2017new}
C.~C. Hammerquist, J.~A. Nairn, A new method for material point method particle
  updates that reduces noise and enhances stability, Computer Methods in
  Applied Mechanics and Engineering 318 (2017) 724--738.

\bibitem{jiang2017angular}
C.~Jiang, C.~Schroeder, J.~Teran, An angular momentum conserving
  affine-particle-in-cell method, Journal of Computational Physics 338 (2017)
  137--164.

\bibitem{holzapfel2002nonlinear}
G.~A. Holzapfel, Nonlinear solid mechanics: A continuum approach for
  engineering science, Meccanica 37~(4) (2002) 489--490.

\bibitem{de2011computational}
E.~A. de~Souza~Neto, D.~Peric, D.~R. Owen, Computational Methods for
  Plasticity: Theory and Applications, John Wiley \& Sons, 2011.

\bibitem{borja2013plasticity}
R.~I. Borja, Plasticity Modeling \& Computation, Springer, 2013.

\bibitem{bardenhagen2004generalized}
S.~G. Bardenhagen, E.~M. Kober, The generalized interpolation material point
  method, Computer Modeling in Engineering and Sciences 5~(6) (2004) 477--496.

\bibitem{harlow1964particle}
F.~H. Harlow, The particle-in-cell computing method for fluid dynamics, Methods
  in Computational Physics 3 (1964) 319--343.

\bibitem{deVaucorbeil2020total}
A.~de~Vaucorbeil, V.~P. Nguyen, C.~R. Hutchinson, {A Total-Lagrangian Material
  Point Method for solid mechanics problems involving large deformations},
  Computer Methods in Applied Mechanics and Engineering 360 (2020) 112783.

\bibitem{deVaucorbeil2021modelling}
A.~de~Vaucorbeil, V.~P. Nguyen, Modelling contacts with a total lagrangian
  material point method, Computer Methods in Applied Mechanics and Engineering
  373 (2021) 113503.

\bibitem{ortiz2015improved}
A.~Ortiz-Bernardin, M.~Puso, N.~Sukumar, {Improved robustness for
  nearly-incompressible large deformation meshfree simulations on Delaunay
  tessellations}, Computer Methods in Applied Mechanics and Engineering 293
  (2015) 348--374.

\bibitem{broccardo2009assumed}
M.~Broccardo, M.~Micheloni, P.~Krysl, Assumed-deformation gradient finite
  elements with nodal integration for nearly incompressible large deformation
  analysis, International Journal for Numerical Methods in Engineering 78~(9)
  (2009) 1113--1134.

\bibitem{cook1974improved}
R.~D. Cook, Improved two-dimensional finite element, Journal of the Structural
  Division 100~(9) (1974) 1851--1863.

\bibitem{wallstedt2008evaluation}
P.~C. Wallstedt, J.~Guilkey, An evaluation of explicit time integration schemes
  for use with the generalized interpolation material point method, Journal of
  Computational Physics 227~(22) (2008) 9628--9642.

\bibitem{hu2019taichi}
Y.~Hu, T.-M. Li, L.~Anderson, J.~Ragan-Kelley, F.~Durand, Taichi: a language
  for high-performance computation on spatially sparse data structures, ACM
  Transactions on Graphics (TOG) 38~(6) (2019) 201.

\bibitem{rodriguez2016particle}
J.~Rodriguez, J.~M. Carbonell, J.~Cante, J.~Oliver, {The particle finite
  element method (PFEM) in thermo-mechanical problems}, International Journal
  for Numerical Methods in Engineering 107~(9) (2016) 733--785.

\bibitem{alkafaji2013formulation}
I.~K. Al-Kafaji, {Formulation of a dynamic material point method (MPM) for
  geomechanical problems}, Ph.D. thesis, University of Stuttgart (2013).

\bibitem{li2020incremental}
M.~Li, Z.~Ferguson, T.~Schneider, T.~R. Langlois, D.~Zorin, D.~Panozzo,
  C.~Jiang, D.~M. Kaufman, Incremental potential contact: intersection-and
  inversion-free, large-deformation dynamics., ACM Transactions on Graphics
  (TOG) 39~(4) (2020) 49.

\bibitem{zhao2022barrier}
Y.~Zhao, J.~Choo, Y.~Jiang, M.~Li, C.~Jiang, K.~Soga, A barrier method for
  frictional contact on embedded interfaces, Computer Methods in Applied
  Mechanics and Engineering 393 (2022) 114820.

\bibitem{li2022bfemp}
X.~Li, Y.~Fang, M.~Li, C.~Jiang, {BFEMP: Interpenetration-free MPM--FEM
  coupling with barrier contact}, Computer Methods in Applied Mechanics and
  Engineering 390 (2022) 114350.

\bibitem{jiang2022hybrid}
Y.~Jiang, Y.~Zhao, C.~E. Choi, J.~Choo, Hybrid continuum--discrete simulation
  of granular impact dynamics, Acta Geotechnica 17 (2022) 5597--5612.

\bibitem{wallstedt2011weighted}
P.~Wallstedt, J.~Guilkey, A weighted least squares particle-in-cell method for
  solid mechanics, International Journal for Numerical Methods in Engineering
  85~(13) (2011) 1687--1704.

\bibitem{kamojjala2015verification}
K.~Kamojjala, R.~Brannon, A.~Sadeghirad, J.~Guilkey, Verification tests in
  solid mechanics, Engineering with Computers 31~(2) (2015) 193--213.

\bibitem{sulsky2016improving}
D.~Sulsky, M.~Gong, Improving the material-point method, Innovative Numerical
  Approaches for Multi-Field and Multi-Scale Problems: In Honor of Michael
  Ortiz's 60th Birthday (2016) 217--240.

\bibitem{zhang2021truncated}
K.~Zhang, S.-L. Shen, A.~Zhou, D.~Balzani, Truncated hierarchical {B}-spline
  material point method for large deformation geotechnical problems, Computers
  and Geotechnics 134 (2021) 104097.

\bibitem{steffen2008examination}
M.~Steffen, P.~Wallstedt, J.~Guilkey, R.~Kirby, M.~Berzins, {Examination and
  analysis of implementation choices within the material point method (MPM)},
  Computer Modeling in Engineering and Sciences 31~(2) (2008) 107--127.

\bibitem{schulz2019consistent}
S.~Schulz, G.~Sutmann, A consistent boundary method for the material point
  method-using imge particles to reduce boundary artefacts, in: PARTICLES VI:
  Proceedings of the VI International Conference on Particle-Based Methods:
  Fundamentals and Applications, CIMNE, 2019, pp. 522--533.

\bibitem{nakamura2023taylor}
K.~Nakamura, S.~Matsumura, T.~Mizutani, Taylor particle-in-cell transfer and
  kernel correction for material point method, Computer Methods in Applied
  Mechanics and Engineering 403 (2023) 115720.

\bibitem{baumgarten2023analysis}
A.~S. Baumgarten, K.~Kamrin, Analysis and mitigation of spatial integration
  errors for the material point method, International Journal for Numerical
  Methods in Engineering 24~(11) (2023) 2449--2497.

\bibitem{bui2021smoothed}
H.~H. Bui, G.~D. Nguyen, {Smoothed particle hydrodynamics (SPH) and its
  applications in geomechanics: From solid fracture to granular behaviour and
  multiphase flows in porous media}, Computers and Geotechnics 138 (2021)
  104315.

\end{thebibliography}

\end{document}